\documentclass[reqno]{amsart}
\usepackage{amsthm,amsmath,amsfonts,amssymb,amscd,mathrsfs,diagrams}
\usepackage{enumerate}
\usepackage{bbm} %\mathbb numbers and other symbols
\usepackage{hyperref}%,supertabular}
\usepackage{eucal}

\newarrow {Dotsto} ....> %makes a dotted arrow 
\newarrow {Onto} ----{>>} %surjective
\newarrow {Into} C---> %injective
\newarrow{Eq} ===={=} %equal
\newcommand{\x}{\mathbf{x}}
\newcommand{\p}{\mathbf{p}}

\newcommand{\Vf}{\mathfrak{C}\kern-.05em\mathfrak{R}}

\renewcommand{\setminus}{\smallsetminus}
\DeclareMathOperator{\Pf}{Pf} %Pfaffian

\DeclareMathOperator{\B}{{B}\kern-.1em}

\DeclareMathOperator{\Fl}{Flag}
\DeclareMathOperator{\GL}{GL}

\DeclareMathOperator{\Gr}{Gr}
\DeclareMathOperator{\Hom}{Hom}

\DeclareMathOperator{\Pic}{Pic}
\DeclareMathOperator{\Proj}{Proj}

\DeclareMathOperator{\SL}{SL}

\DeclareMathOperator{\Spec}{Spec}

\DeclareMathOperator{\crit}{Crit}

\DeclareMathOperator{\pt}{pt}

\DeclareMathOperator{\stab}{Stab}

\newcommand{\0}{\boldsymbol{0}}
\newcommand{\Acal}{{\mathscr A}}
\newcommand{\Bcal}{{\mathscr B}}
\newcommand{\Bchi}{\B\kern-.1em\chi}
\newcommand{\BchiR}{\B\kern-.1em\chiR}

\newcommand{\CC}{\mathbb C}
\newcommand{\Ccal}{\mathscr C}
\newcommand{\Ecal}{\mathscr{E}}

\newcommand{\Gammah}{\widehat{\Gamma}} %\character group of \Gamma
\newcommand{\Gh}{\widehat{G}} %character group of G
\newcommand{\GQ}{B} %matrix of G weights

\newcommand{\Hcal}{\mathscr H}

\newcommand{\LGQ}{{\operatorname{LGQ}}} %LG quasimaps (without P-fields)
\newcommand{\Lcal}{\mathscr L}

\newcommand{\Ocal}{\mathscr O}
\newcommand{\PP}{\mathbb P}
\newcommand{\Pcal}{\mathscr P}

\newcommand{\QQ}{\mathbb Q}

\newcommand{\Rq}{c} % the R-charge weights

\newcommand{\LV}{\mathbf{L}} %line bundle on V defined by character \theta of G

\newcommand{\WP}{W\PP} %Weighted projective space

\newcommand{\XLB}{{L}}  %Line bundle on X_\theta defined by a character of G

\newcommand{\XS}{\X^{\operatorname{sympl}}}
\newcommand{\Xcal}{\mathscr{X}}

\newcommand{\X}{\Xcal}

\newcommand{\ZZ}{\mathbb Z}

\newcommand{\aff}{/_{\!\!\operatorname{aff}}}
\newcommand{\age}{\operatorname{age}}

\newcommand{\Gdeg}{\deg}%^{\Gamma}}  %degree map on chars of Gamma

\newcommand{\chat}{\hat{c}}
\newcommand{\chiR}{\zeta}

\newcommand{\cjcl}[1]{[ \kern-.15em [ #1 ]\kern-.15em ]}

\newcommand{\crst}{\mathscr{C \kern-.25em R}} %stack quot of critical locus in \X

\newcommand{\dsand}{{\quad \text{and} \quad}} %put an "and" with space on a displayed line.

\newcommand{\genj}{\langle J \rangle} %the group generated by J

\newcommand{\git}[1]{{\!/\!\!/_{\kern-.2em #1}}} %git quotient
\newcommand{\sympl}[1]{{\!/\!\!/_{\kern-.2em #1}^{\operatorname{spl}}}} %symplectic reduction

\newcommand{\hfrak}{\mathfrak{h}}

\newcommand{\kk}{\omega} % the canonical bundle
\newcommand{\pkk}{\mathring{\omega}}  %principal bundle assoc to omega
\newcommand{\klog}{\kk_{\log}}
\newcommand{\klogc}{\kk_{\log,\Ccal}}    %the log canonical bundle on C
\newcommand{\pklogc}{\pkk_{\log,\Ccal}}    %principle bundle assoc to log canonical bundle on C

\newcommand{\lift}{\vartheta} %lift of character on G to Gamma
\newcommand{\mmapvalue}{\tau}
\newcommand{\mmap}{\mu}
\newcommand{\mrkp}{y} %marked points.

\newcommand{\qmp}{\mathcal{Q}} %an LG quasimap

\newcommand{\spl}{\xi} % map \Gamma \to G/J
\newcommand{\spn}{\varkappa}  %isom of \zeta_*\Pcal to \pklogc

\newcommand{\thetaweight}{e} % -weight of the characater \theta

\newcommand{\trp}{^{\mathsf T}}  %matrix transpose

\renewcommand{\u}{\sigma} %the section defining the LGquasimap
\newcommand{\ufrak}{\mathfrak{u}}

\newcommand{\ve}{\varepsilon}

\renewcommand{\Re}{\mathfrak{Re}}

\renewcommand{\hom}{\operatorname{Hom}}

\newcommand{\Q}{\QQ}

\newcommand{\C}{\CC}

\newtheorem{thm}{Theorem}[subsection] % number like 3.1, 3.2, 3.3, etc.
\newtheorem{lem}[thm]{Lemma} % numbered with thm
\theoremstyle{definition} % 'here we change the style'
\newtheorem{defn}[thm]{Definition} % numbered with thm

\theoremstyle{remark} % 'style changed again'
\newtheorem{rem}[thm]{Remark}% numbered with thm
\begin{document}

\title{The Moduli Space in the Gauged Linear Sigma Model}
\author{Huijun Fan, Tyler Jarvis and Yongbin Ruan}
\date{\today}
\thanks{The second author's work was partially supported by NSF grant DMS-1564502. The third author's work was partially supported by NSF grants DMS 1159265
and DMS 1405245.}
\begin{abstract}
This is a survey article for the mathematical  theory of Witten's Gauged Linear Sigma Model, as developed recently by the authors.
Instead of developing the theory in the most general setting, in this paper we focus on the description of the moduli.
\end{abstract}

\maketitle
\setcounter{tocdepth}{1}
\tableofcontents

\section{Introduction}

In 1991, in an effort to generalize his famous conjecture regarding the KdV-hierarchy for
the intersection theory of the moduli space of Riemann surfaces \cite{Wit:92a, Wit:93}, Witten 
proposed a remarkable first-order, nonlinear, elliptic PDE associated to an
arbitrary quasihomogeneous singularity. It has the simple form
\begin{equation}\label{eq:Witten}
\bar{\partial}u_i+\overline{\frac{\partial W}{\partial u_i}}=0,
\end{equation}
where $W$ is a quasihomogeneous polynomial with an isolated singularity at the origin, and $u_i$ is interpreted as a section of an appropriate orbifold line bundle on an orbifold Riemann
surface $\Ccal$.
  
During the last decade we have carried out a comprehensive treatment of the Witten
equation and have used it to construct a theory similar to Gromov-Witten theory \cite{FJR:07b, FJR:08, FJR:07a}.  This so-called FJRW-theory  can be viewed as the
Landau-Ginzburg phase of a Calabi-Yau hypersurface
\[
X_W=\{W=0\} \subset \WP^{n-1}
\]
in weighted projective space.  The
relation between the Gromov-Witten theory of $X_W$ and the FJRW-theory
of $W$ is the subject of the Landau-Ginzburg/Calabi-Yau
correspondence, a famous duality from physics. More recently, the
LG/CY correspondence has been reformulated as a precise
mathematical conjecture \cite{Rua:12}, and a great deal of  progress has
been made on this conjecture \cite{CIR:12, ChiRu:10, ChiRu:11, PS, LPS}.

A natural question is whether the LG/CY correspondence can be
generalized to complete intersections in projective space, or more
generally to toric varieties. The physicists' answer is ``yes.'' In
fact, Witten considered this question in the early 90s \cite{Wit:92a}
in his effort to give a physical derivation of the LG/CY
correspondence. In the process, he invented an important model in
physics called the \emph{Gauged Linear Sigma Model (GLSM)}. From the
point of view of partial differential equations, the gauged linear
sigma model generalizes the Witten Equation \eqref{eq:Witten} to the
\emph{Gauged Witten Equation}
\begin{align}\label{eq:GaugedWitten}
\bar{\partial}_A u_i+\overline{\frac{\partial W}{\partial
    u_i}}&=0,\\ *F_A&=\mmap,
\end{align}
where $A$ is a connection of certain principal bundle, and $\mmap$ is the moment map of the
GIT-quotient, viewed as a symplectic quotient.  In general, both the
Gromov-Witten theory of a Calabi-Yau complete intersection $X$ and the
LG dual of $X$ can be expressed as gauged linear sigma
models. Furthermore, the LG/CY correspondence can be interpreted as a
variation of the moment map $\mmap$ (or a deformation of GIT) in the GLSM. 

During last several years,
we constructed a rigorous 
mathematical theory for the gauged linear sigma model \cite{FJR:15a,FJR:15b}, and this new model seems to have many applications (see, for example \cite{RR,RRS,CJR}).  Our new theory is a generalization of FJRW-theory
from the case of a finite gauge group to the case where the gauge group is any reductive Lie group. Surprisingly, our older theory (finite group) is one of the more difficult cases in the general theory. 

We deal with the Gauged Witten Equation both analytically \cite{FJR:15b} and algebraically \cite{FJR:15a}.  But in this paper we focus on the algebraic version of the theory.
In the case of a continuous group, we can treat the Gauged Witten Equation algebraically using some stability conditions. It turns out to be very convenient to incorporate the stability conditions from the quasimap theory of Ciocan-Fontanine, Kim, Maulik, and Cheong \cite{CCFK:14,CFKM:11,CFKi:10, Kim:11}.  And indeed, it is natural to view our GLSM-theory as a union of FJRW-theory with quasimap theory.  However, it is possible to impose other stability conditions (see \cite{CLL:15, ChKi:15}).

Understanding the details of the mathematical construction of the GLSM can be daunting. This paper is an attempt to help the reader navigate past the more technical constructions and begin to understand the underlying moduli spaces.  In the next section we give a brief overview of the main ingredients to the theory, and in the subsequent section we briefly review the definition of the moduli problem and the most important stability conditions.  The remainder of the paper is focused on giving examples of the moduli spaces for various choices of input data.

\section{The Basic Setting}

The \emph{input data} of our new theory consists of the following.  We discuss these in more detail below. 
\begin{enumerate}
\item A finite dimensional vector space $V$ over $\CC$.
\item A reductive algebraic group $G\subseteq GL(V)$, sometimes called the \emph{gauge group}.
\item A $G$-character $\theta$ with the property $V^s_G(\theta)=V^{ss}_G(\theta)$. We say that it defines a \emph{strongly regular phase} $\X_\theta = [V\git{\theta}G]$. 
\item A choice of $\CC^*$ action on $V$, called the \emph{$R$-charge} and denoted $\CC_R^*$.  This action is required to commute with the $G$-action, and we require $G\cap \CC^*_R = \langle J \rangle$ to have finite order $d$.  
\item A $G$-invariant quasihomogeneous polynomial $W:V\to \CC$, called the \emph{superpotential}, having degree $d$ with respect to the $\CC^*_R$ action. We require that the GIT quotient $\crit(W)\git{\theta} G$ of the critical locus $\crit(W)$ be compact.
\item A stability parameter $\ve$.  This can be any positive rational number, but in practice, the two most useful cases are the limiting cases of $\ve \to \infty$ or $\ve\to 0+$.  Fortunately, these cases are also easier to describe than the positive rational cases.  For simplicity, in this paper we will only discuss the limiting cases of $0+$ and $\infty$.

\item A $\Gamma$-character $\lift$, where $\Gamma$ is the subgroup of $\GL(V)$ generated by $G$ and $\CC^*_R$ by $\Gamma$.  We require that $\lift$ define a \emph{lift} of $\theta$, meaning that $\lift|_G= \theta$.  Except in the case of $\ve=0^+$, we also require that this lift be a \emph{good lift}, meaning that it satisfy $V^{ss}_{\Gamma}(\lift)=V^{ss}_G(\theta)$.  A choice of good lift affects the stability conditions for the moduli space.  But in the case of $\ve=0+$ the lift need not be good, and every lift will produce the same stability conditions and the same moduli space.
\end{enumerate}
With the above input data we construct a theory with the following main ingredients: 

\begin{enumerate}
\item \emph{A state space}, which is the relative Chen-Ruan cohomology of the quotient $\X_\theta = [V\git{\theta} G]$ with an
additional shift by $2q$.  For each conjugacy class $\Psi\subset G$, let 
\[
I(\Psi) = \{(v,g)\in V^{ss}_\theta \times G |  g \in \Psi\}
\]
and 
\[
\X_{\theta,\Psi} = [I(\Psi)/G].
\]
The state space is  
\[
\Hcal_{W, G}=\bigoplus_{\alpha\in \QQ}\Hcal^{\alpha}_{W, G} = \bigoplus_{\Psi}\Hcal_{\Psi},
\]
where the sum runs over those conjugacy classes $\Psi$ of $G$ for which  $\X_{\theta,\Psi}$ is nonempty, and where 
\[
\Hcal^{\alpha}_{W, G}=H^{\alpha+2q}_{CR}(\X_{\theta},W^{\infty},
 \QQ)=\bigoplus_{{\Psi}} H^{\alpha-2\age{(\gamma)}+2q}(\X_{\theta,\Psi},
  W^{\infty}_{\Psi},\QQ),
\]
  and 
\[
\Hcal_{\Psi}=H^{\bullet+2q}_{CR}(\X_{\theta,\Psi},W^{\infty},
 \QQ)= \bigoplus_{\alpha\in \QQ} H^{\alpha-2\age{(\gamma)}+2q}(\X_{\theta,\Psi},
  W^{\infty}_{\Psi},\QQ).
\] 
Here $W^\infty = \Re(W)^{-1}(M,\infty) \subset [V\git{\theta} G]$ for some large, real $M$. 

\item 
\emph{The stack of LG-quasimaps:}

We denote by $\crst_\theta = [\crit_G^{ss}(\theta)/G] \subset  [V\!\git{\theta} G] = [V^{ss}_G(\theta)/G]$ 
the GIT quotient (with polarization $\theta$) of the critical locus of $W$. Our main object of study is the stack 
\[
\LGQ_{g,k}^{\ve,\lift}(\crst_{\theta}, \beta)
\] 
of \emph{$(\ve,\lift)$-stable Landau-Ginzburg quasimaps to $\crst_\theta$}.  

\item \emph{A virtual cycle:}
\[
[\LGQ_{g,k}^{\ve,\lift}(\crst_\theta, \beta)]^{vir}\in H_*(\LGQ_{g,k}^{\ve,\lift}(\crst_\theta, \beta), \QQ)
\]
with virtual dimension
\[
\dim_{vir}=\int_{\beta} c_1(V\git{\theta} G)+(\chat_{W,G}-3)(1-g)+k-\sum_i (\age(\gamma_i)-q),
\]
where $\chat_{W, G}$ is the \emph{central charge} (see Definition~\ref{def:chat}).

\item \emph{Numerical invariants:}
Using the virtual cycle, we can define correlators  
\[
\langle \tau_{l_1}(\alpha_1),\cdots, \tau_{l_k}(\alpha_k)\rangle=\int_{[\LGQ_{g,k}^{\ve,\lift}(\crst_\theta,\beta)]^{vir}}\prod_i ev_i^*(\alpha_i)\psi^{l_i}_i.
\] 
One can then define a generating function in the standard fashion.
These invariants satisfy the usual gluing axioms whenever all insertions are of compact type.
\end{enumerate}

In the rest of this section we will discuss some of the input data and the state space in more detail.  

\subsection{GIT and Symplectic Quotients}
The first two pieces of data consist of a reductive algebraic group $G$ (the \emph{gauge group}) acting on a finite-dimensional vector space $V \cong \CC^{n} $.  We do not require $G$ to be connected, but we require that  $G/G_0$ be finite, where  $G_0$ is the identity component of $G$.  If the gauge group action on $V$ factors through $\SL(V)$ then we say that it satisfies the \emph{Calabi-Yau condition.}  But in general we do not require that $G$ satisfy this condition.

We wish to consider the quotient stack $[Z/G]$ for a closed subvariety $Z\subseteq V$, but since the group $G$ may not be compact, the quotient is not generally separated (Hausdorff). Geometric Invariant Theory (GIT) and symplectic reduction each give a way to construct separated quotients.

\subsubsection{GIT Quotients} 
The key to constructing a separated quotient using GIT is to choose a linearization of the action of $G$ on $Z$, i.e., a lifting of the action of $G$ to a line bundle $\LV$ over $Z$.  We always assume that the linearization on $Z$ is induced by a linearization on $V$.  Since $V=\CC^n$, any line bundle $\LV$ on $V$ is trivial $\LV =
V\times \CC$, and the linearization is determined by a character
$\theta:G\to \CC^*$.
\begin{defn}
For any character $\theta:G\to \CC^*$ we write $\LV_\theta$ for the
line bundle $V\times \CC$ with the induced linearization.  We also often write $\LV_\theta$ to denote the corresponding line bundle on $Z$.
\end{defn}

Geometric Invariant Theory identifies an open subset $Z^{ss}(\theta)$ of
\emph{$\LV_\theta$-semistable} points in $Z$ as the set of those
points $v\in Z$ for which there exists a positive integer $k$ and a $G$-invariant section $f\in H^0(Z,\LV_\theta^{\otimes k})^G$ such that $f(v)\neq
0$.  We denote the set of points in $Z$ that are semistable with respect to $G$ and $\theta$ by $Z^{ss}_G(\theta)$.
The GIT quotient stack $[Z\git{\theta} G]$ is defined
to be the stack
\[
[Z\git{\theta} G] = [Z^{ss}_G(\theta)/G].
\]
Let $Z\aff G$ be the \emph{affine quotient} given by $Z\aff G = \Spec(\CC[Z^*]^G)$, where $\CC[Z^*]$ is the ring of regular functions on $Z$.
The GIT quotient stack $[Z\git{\theta} G]$
is an algebraic stack with an underlying
(coarse moduli) space 
\[
Z\git{\theta} G = Z^{ss}_G(\theta)/G = 
\Proj_{Z\aff G}\left(\bigoplus_{k\ge 0} H^0(Z,\LV_\theta^k)^G\right).
\]
The linearization $\LV_{\theta}$ induces a line bundle (a.k.a. a \emph{polarization}) on
$[Z\git{\theta} G]$, which we denote by $\XLB_{\theta}$.  

\begin{defn}
We say that a point $v\in V$ is \emph{stable} with respect to the linearization $\theta$ (or $\theta$-stable) if 
\begin{enumerate}
\item $v$ is $\theta$-semistable
\item The stabilizer $\stab_G(v) = \{g\in G\mid gv = v\}$ is finite.
\end{enumerate}
We denote the set of $\theta$-stable points of $Z$ by $Z^s_G(\theta)$.  We say that a point is \emph{unstable} if it is not semistable.
\end{defn}
The stable locus is important because the quotient stack $[Z^{s}_G(\theta)/G]$ is a Deligne-Mumford stack, whereas  $[Z^{ss}_G(\theta)/G]$ is not necessarily Deligne-Mumford.

\begin{rem}
Mumford-Fogarty-Kirwan\cite{MFK:94} use the name \emph{properly stable} to describe what we call stable.
\end{rem}
\begin{rem}
For any integers $\ell,k>0$, each $f\in H^0(Z,\LV_\theta^{\otimes k})^G$ also satisfies $f^\ell\in H^0(Z,\LV_\theta^{\otimes k \ell})^G$, so it makes sense to extend the GIT constructions to fractional linearizations, corresponding to fractional characters in $\Gh_\Q = \hom(G,\CC^*)\otimes_\ZZ \QQ$.
\end{rem}

For a fixed $Z$, changing the linearization gives a different
quotient.  The space of (fractional) linearizations is divided into chambers, and
any two linearizations lying in the same chamber have isomorphic GIT
quotients.  We call the isomorphism classes of these quotients
\emph{phases}.  If the linearizations lie in distinct chambers, the
quotients are birational to each other, and are related by flips
\cite{Tha:96, DoHu:98}.

\begin{defn}\label{def:StronglyRegular}
We say that $\theta\in\Gh_\QQ$ (or the corresponding linearization $\LV_\theta$) is \emph{strongly regular} 
if $V^{ss}_G(\theta)$ is not empty and $V^{s}_G(\theta) = V^{ss}_G(\theta)$.
\end{defn}
For purposes of this paper, all linearizations need to be \emph{strongly regular}.

\begin{rem}
For any strongly regular phase $\theta$, the complex dimension of $\X_\theta = [V\!\git{\theta} G]$ is $n-\dim(G)$.
\end{rem}

\subsubsection{Symplectic Quotients}

One may also think of the GIT quotients as symplectic
reductions.  Take $Z \subseteq \CC^{{n}}$ with the 
standard K\"ahler form $\omega=\sum_i dz_i\wedge
d\bar{z}_i.$ Since $G$ is reductive, it is the complexification of a
maximal compact Lie subgroup $H$, acting on $Z$ via a faithful unitary
representation $H\subseteq U({{n}})$. Denote the Lie algebra of $H$ by ${\hfrak}$.

We have a Hamiltonian action of $H$ on $Z$ with moment map $\mmap_Z: Z
\to {\hfrak}^*$ for the action of $H$ on $Z$, given by
\[
\mmap_Z(v)(Y) = \frac{1}{2} \overline{v}\trp Y v = \frac{1}{2}
\sum_{i,j \le {{n}}} \bar{v}_i Y_{i,j} v_j
\]
for $v \in Z$ and $Y\in {\hfrak}$.  If $\mmapvalue \in {\hfrak}^*$ is a
 value of the moment map, then the locus $\mmap^{-1}(H\mmapvalue)$ is
an $H$-invariant set, and the symplectic orbifold quotient of $Z$ at
$\mmapvalue$ is defined as
\[
\left[Z\sympl{\mmapvalue} H\right] = \left[ \mmap_Z^{-1}(H\mmapvalue)/H
  \right] = \left[ \mmap_Z^{-1}(\mmapvalue)/H_\mmapvalue
  \right],
\]
where $H_\mmapvalue$ is the stabilizer in $H$ of $\mmapvalue$.
The value $\mmapvalue$ of the moment map plays the role for the symplectic quotient that the linearization $\theta$ plays for the GIT quotient.  These are related by the following result.

\begin{thm}[{\cite[Cor.~2.1.8]{FJR:15a}}]
Let $\theta\in\Gh$ be a character of $G$.  Taking
derivations of the character $\theta$ defines a weight
$\mmapvalue_\theta\in\hfrak^*$.  Whenever the coadjoint orbit of $\tau_\theta$ in $\hfrak^*$ is trivial (e.g., in the case that $\tau$ is in the Lie algebra of the center of $G$, or if $G$ is Abelian), then  we have 
\[
\left[Z\sympl{-\mmapvalue_\theta} H\right]
= [Z\git\theta G ].
\]
\end{thm}

As with the space of GIT linearizations, the space $\mathfrak{h}^*$  is divided into chambers; and values of  $\mmapvalue$  that lie in the same chamber define isomorphic quotients.  The walls between the chambers correspond to the critical points of the moment map $\mu$. In many cases these are easier to identify in the symplectic formulation than in the GIT formulation.

\subsection{Superpotential and Critical Locus}

The next piece of data required for the GLSM is the \emph{superpotiential}, which is $G$-invariant polynomial $W:V\to \C$.  We are especially interested in the critical locus of the superpotential.

\begin{defn}\label{defn:WGnondegen}
Let $\theta:G\to \CC^*$ define a  strongly regular phase $\X_\theta = \left[V\!\git{\theta} G\right]$. 
The superpotential $W$ descends to a holomorphic function $W{\colon} \X_{\theta}\rightarrow
\CC$.  Let  $\crit^{ss}_G(\theta) = \{v\in V^{ss}_G(\theta) \mid \frac{\partial W}{\partial x_i} \text{for
  all $i =1, \dots, {{n} }$}\} \subset V^{ss}$ denote the semistable points of the critical
locus.  The group $G$ acts on $\crit^{ss}_G(\theta)$, and the stack quotient is
\[
\crst_\theta = 
\left[\crit^{ss}_G(\theta)/G\right] = \{x\in \X_{\theta} \mid d W = 0 \}\subset \X_{\theta},
\]
where $d W\colon T\X_{\theta} \to T\CC^*$ is the differential of $W$ on $\X_\theta$.  
We say that the pair $(W, G)$ is \emph{nondegenerate
  for $\X_{\theta}$} if the critical locus $\crst_\theta \subset \X_{\theta}$  is
  compact.
\end{defn}

\subsection{$R$-charge and the Group $\Gamma$}\label{sec:GLSM}

The Gauged Linear Sigma Model (GLSM) requires an additional $\CC^*$-action on $V$ 
called the \emph{R-charge}. 
The $R$-charge is a $\CC^*$-action on $V$ of the form $(z_1, \cdots, z_{{n}})
\mapsto (\lambda^{c_1}z_1, \cdots, \lambda^{c_{{n} }}z_{{n} })$.  We
denote this action by $\CC^*_R$  in order to distinguish it from other $\CC^*$ actions (for example when $G=\CC^*)$.  We think of $\CC_R^*$ as a subgroup of
$\GL(V,\CC)$.  This means we require $\gcd(c_1,\dots,c_{n} ) = 1$.  Unlike the case of FJRW theory, we allow the weights $c_i$ of $\CC^*_R$ to be zero or negative.

\begin{rem}
Our choice of $\CC^*_R$-action differs from what what the physics literature calls
\emph{R-charge} by a factor of $2$.  More precisely, the physicists' R-charge is the $\CC^*$-action given by the weights $(2c_1/d,\dots,2c_{{n} }/d)$; but for our purposes, $\CC^*_R$ is the more natural choice.
\end{rem}

 \begin{defn}
We define the \emph{exponential grading element} $J\in \CC^*_R$ to be 
\begin{equation}\label{eq:J}
J = (\exp(2\pi i c_1/d), \dots, \exp(2 \pi i c_{{n} }/d)),
\end{equation}
which has order $d$.
 
It is sometimes convenient to write $q_i = c_i/d$ and $q = \sum_{i=1}^n q_i$ so that 
\[
J = (\exp(2\pi i q_1), \dots, \exp(2 \pi i q_{{n} })).
\]
\end{defn} 
 
We require the actions of $G$ and $\CC^*_R$ to be \emph{compatible}, by which we mean
\begin{enumerate}
\item They commute: $g r = r g$ for any $g\in G$ and any $r \in \CC^*_R$.
\item We have $G \cap \CC^*_R=\genj$.
\end{enumerate}

\begin{defn}
The group $\Gamma$ is the subgroup of $\GL(V)$ generated by $G$
and $\CC_R^*$.
\end{defn}

If $G$ and $\CC^*_R$ are compatible, then every element $\gamma$ of
$\Gamma$ can be written as $\gamma = g r$ for $g\in G, r\in \CC^*_R$; that is,
\[
\Gamma = G \CC_R^*.
\]

The representation $\gamma = gr$ is unique up to an element of $\genj
$. Moreover, there is a well-defined homomorphism 
\begin{align}\label{eq:defzeta}
\chiR{\colon} \Gamma = G\CC^*_R &\rightarrow \CC^*\\
 g(\lambda^{c_1}, \cdots, \lambda^{c_{{n} }}) &\mapsto \lambda^d. \notag
\end{align}
We denote the target of $\chiR$ by $H = \chiR(\CC^*_R)= \CC^*$, to distinguish it from $\CC^*_R$. 
This gives the following
exact sequence:
\begin{equation}\label{eq:chiR-seq}
1 \rTo G \rTo \Gamma \rTo^{\chiR}H \rTo 1
\end{equation}
Moreover, there is another homomorphism 
\begin{align}\label{eq:defspl}
\spl{\colon}\Gamma &\rightarrow  G/\genj  \\
g r &\mapsto g\genj.\notag
\end{align}
This is also well defined, and gives another exact
sequence:
\[
1 \rTo \CC^*_R \rTo \Gamma \rTo^{\spl}G/\genj \rTo 1.
\]

\begin{defn}\label{def:chat}
Let $N = n-\dim(G)$.  We define the \emph{central charge} of the theory for the choices $V,G,\CC^*_R, W$ to be 
\begin{equation}\label{eq:def-chat}
\chat_{W,G}=N-2\sum_{j=1}^{n}  c_j/d = N-2q. 
\end{equation}
\end{defn}

\subsection{Lifts of the Linearization to $\Gamma$}\label{subsec:lifts}

Although we are primarily interested in the GIT quotients of $V$ by $G$, our constructions also depend heavily on the GIT quotients of $V$ by $\Gamma$.  For this, we need a \emph{lift} of the $G$-linearization to a $\Gamma$-linearization.  That is, we require a $\Gamma$-character $\lift$ that {lifts} $\theta$, meaning that $\lift|_G= \theta$.  It is not hard to prove that lifts always exist, including the trivial lift $\lift(gr) = \theta(g)$.

For a given lift $\lift$, we always have $V^{ss}_\Gamma(\lift) \subset V^{ss}_G(\theta)$, but equality does not necessarily hold.  If it does hold, we say that $\lift$ is a \emph{good lift} of $\theta$.  
For the stability parameter $\ve=\infty$ we require a the lift to be a good lift.  For the choice $\ve=0+$, the lift need not be good, and every lift will produce the same theory.
Unfortunately, not every $\theta\in \Gh$ has a good lift for every choice of ($G$-compatible) $\CC^*_R$-action, but most of the interesting examples of GLSMs have a good lift.

\subsection{Choice of $\CC^*_R$}

All of our constructions ultimately depend not on $\CC^*_R$, but rather only on the embeddings $G \subseteq \Gamma\subseteq GL(V)$, on the sum $q = \sum_{i=1}^{n}  q_i = \sum_{i=1}^{n}  c_i/d$ of the $\CC^*_R$ weights, and on a choice of a lift $\lift:\Gamma \to \CC^*$ of $\theta$.

Of course the choice of $q$ and the embedding of $\Gamma$ in $\GL(V)$ put many constraints on $\CC^*_R$; but they still allow some flexibility.
For example, when the gauge group $G$ is a torus with a Calabi-Yau weight system (that is, if its weight matrix $\GQ=(b_{i j})$ satisfies $\sum_j b_{ij} = 0$ for each $i$), then we have a lot of flexibility.  The following lemma is not hard to prove (see \cite[Lem 3.3.1]{FJR:15a}). 
\begin{lem}\label{lem:CstarR}
If the gauge group $G$ is a torus with weight matrix $\GQ=(b_{i j})$, and if we have a compatible $\CC^*_R$ action with weights $(\Rq_1,\dots,\Rq_{n} )$, such that $W$ has $\CC^*_R$-weight $d$, then for any $\QQ$-linear combination $(b'_1, \dots, b'_{n} )$ of rows 
of the gauge weight matrix $\GQ$, we define a new choice of R-weights $({\Rq'_1},\dots,{\Rq'_{n} }) = (\Rq_1+b'_1, \dots, \Rq_{n} +b'_{n} )$.  Denote the corresponding $\CC^*$ action by $\CC^*_{R'}$. 

Since the group $\Gamma'$ generated by $G$ and $\CC^*_{R'}$ lies inside the maximal torus of $\GL({n} ,\CC)$, it is Abelian; and so we automatically have that $G$ and $\CC^*_{R'}$ commute.  We also have the following:
\begin{enumerate}
\item The group $\Gamma'$ generated by $G$ and $\CC^*_R$ is the same as the group $\Gamma$ generated by $G$ and $\CC^*_R$.
\item The $\CC^*_{R'}$-weight of $W$ is equal to $d$.
\item $G \cap \CC^*_{R'} = G\cap \CC^*_R = \genj$, where $J$ is the element defined by Equation~\eqref{eq:J} for the original $\CC^*_R$ action.
\item If $\GQ$ is a Calabi-Yau weight system, then  for both $\CC^*_R$ and $\CC^*_{R'}$ the sum of the weights $q = \sum q_i = \sum c_i/d$ is the same and the central charge $\chat_W$ is the same.
\end{enumerate}
\end{lem}

\subsection{Hybrid Models}

A very important subclass of the toric examples---when $G = (\C^*)^m$---consists of the so-called \emph{hybrid models}.  Several examples of the hybrid model have been worked out in detail by E.~Clader in \cite{Cla:13}.

\begin{defn}
For a torus $G=(\CC^*)^m$, a phase $\theta$ of $(W, G)$ is called a \emph{hybrid model} if 
\begin{enumerate}
\item The quotient $\X_{\theta}\rightarrow \X_{base}$ has the structure of a toric bundle over a compact base $\X_{\mathrm{base}}$, and 
\item The $\CC^*_R$-weights of the base variables are all zero. 
\end{enumerate}
\end{defn}

\section{The GLSM Moduli Space}

Given a choice of $V$, $G$, $W$, $\CC^*_R$, $\theta$, and an additional choice of stability parameters $\ve,\lift$, the ``moduli space'' for the GLSM is the stack of \emph{$(\ve,\lift)$-stable Landau-Ginzburg quasimaps} to the critical locus $\crst_\theta$ of $W$, which we describe below.  This space is naturally a substack of the stack of \emph{$(\ve,\lift)$-stable Landau-Ginzburg quasimaps} to $\X_\theta$, and that larger space plays an important role in the construction of the virtual class---similar to the role of $p$-fields for Gromov-Witten theory.

\subsection{Landau-Ginzburg Quasimaps}

\begin{defn}
For any $k$-pointed orbicurve $\Ccal,\mrkp_1,\dots, \mrkp_k$, denote by $\klogc$ the line bundle 
\[
\klogc = \omega_{\Ccal} \otimes \Ocal(\mrkp_1) \otimes \cdots \otimes \Ocal(\mrkp_k) = \Ocal\left(K_\Ccal + \sum_{i=1}^k \mrkp_i\right),
\]
where $\omega_\Ccal = \Ocal(K_\Ccal)$ is the canonical bundle on $\Ccal$.

Also, let $\pklogc$ denote the principal $\C^*$-bundle on $\Ccal$ corresponding to the line bundle $\klogc$. 
\end{defn}

\begin{defn}\label{defn:LGquasimap}
A \emph{prestable, $k$-pointed, genus-$g$,
 LG-quasimap to $\X_\theta$} is a tuple $(\Ccal,\mrkp_1,\dots,\mrkp_k, \Pcal, u, \spn)$ consisting of
\begin{enumerate}[     A.) ]
   \item A prestable, $k$-pointed orbicurve $(\Ccal, \mrkp_1,\dots,\mrkp_k)$
     of genus $g$.
   \item A principal (orbifold) $\Gamma$-bundle $\Pcal{\colon} \Ccal \to
     \B\Gamma$ over $\Ccal$.
   \item A global section $\u{\colon} \Ccal \to \Ecal = \Pcal\times_{\Gamma} V$.
   \item An isomorphism $\spn \colon \chiR_*\Pcal \to \pklogc$ of principal $\CC^*$-bundles.
   % ($\pklogc$ indicates the principle $\CC^*$-bundle associated to the line bundle $\klogc$).
   
\end{enumerate}
 such that
\begin{enumerate}
    \item The morphism of stacks $\Pcal{\colon}\Ccal \to \B\Gamma$ is
      representable (i.e., for each point $\mrkp$ of $\Ccal$, the induced
      map from the local group $G_\mrkp$ to $\Gamma$ is injective).
    \item \label{item:basepoints}
      The set $B$ of points $b \in \Ccal$ such that any point $p$ of the fiber $\Pcal_b$ over $b$ is mapped by $\u$
      into an $\LV_\theta$-unstable $G$-orbit of $V$ is finite, and this set is 
      disjoint from the nodes and marked points of $\Ccal$.
\end{enumerate}

A \emph{prestable, $k$-pointed, genus-$g$,
 LG-quasimap to $\crst_\theta$} is a  \emph{prestable, $k$-pointed, genus-$g$, LG-quasimap to $\X_\theta$} such that the image  of the induced map $[\u]{\colon}\Pcal \to V$ lies 
    in $\crit(W)$.   
\end{defn}

\begin{defn}
The points $b$ occurring in condition (\ref{item:basepoints}) above are called 
\emph{base points} of the quasimap.  That is, 
$b\in\Ccal$ is a base point if at least one point (and hence every point) of the fiber $\Pcal_b$ over $b$ is mapped by $\u$
      into an $\LV_\theta$-unstable $G$-orbit of $V$.
\end{defn}

\begin{defn}\label{def:CLBchi}
For any prestable LG-quasimap $\qmp = (\Ccal, \mrkp_1,\dots,\mrkp_k, \Pcal, \u, \spn)$, a $\Gamma$-equivariant line bundle $\LV \in \Pic^{\Gamma}(V)$ determines a line bundle $\Lcal = \Pcal\times_\Gamma \LV$ over $\Ecal = \Pcal\times_\Gamma V$, and pulling back along $\u$ gives a line bundle $\u^*(\Lcal)$ on $\Ccal$.  

In particular, any character $\alpha\in\Gammah = \Hom(\Gamma,\C^*)$ determines a $\Gamma$-equivariant line bundle $\LV_\alpha$ on $V$ and hence a line bundle $\u^*(\Lcal_\alpha)$ on $\Ccal$.  
\end{defn}

\subsection{Stability Conditions for the Stack of LG-Quasimaps}

\begin{defn}\label{def:LGQ-deg}
For any $\alpha \in \Gammah$ define the \emph{degree} of $\alpha$ on $\qmp$ to be  
\[
\Gdeg_\qmp(\alpha) = {\deg_{\Ccal}(\u^*(\Lcal_{\alpha}))} \in \QQ.
\]
This defines a homomorphism $\Gdeg_\qmp:\Gammah \to \QQ$.

For any $\beta\in \Hom(\Gammah,\QQ)$ we say that an 
LG-quasimap  $\qmp = (\Ccal,x_1,\dots,x_k,\Pcal,\u,\spn)$ has \emph{degree  $\beta$} if $\Gdeg_\qmp = \beta$.
\end{defn}

\begin{defn}
Given a polarization $\theta\in \Gh$, a lift $\lift\in \Gammah$ of $\theta$ and a  prestable LG-quasimap $\qmp = (\Ccal,x_1,\dots,x_k,\Pcal,\u,\spn)$, we say that $\qmp$ is \emph{$0+$-stable} if 
\begin{enumerate}
\item Every rational component has at least two special points (a mark $\mrkp_i$ or a node), and 
\item\label{it:zero-plus2} On every irreducible component $\Ccal'$ with trivial $\omega_{log,\Ccal'}$, the line bundle $\u^*(\Lcal_{\lift})$ has positive degree.
\end{enumerate}
\end{defn}
It turns out that condition~(\ref{it:zero-plus2}) holds for one lift if and only if it holds for all lifts, because for any two lifts $\lift$ and $\lift'$ of $\theta$, the bundles $\u^*(\Lcal_\lift)$ and $\u^*(\Lcal_{\lift'})$ always differ by a power of $\klogc$ (see \cite[Prop.~4.2.14]{FJR:15a}).  Moreover, over $\Q$ the ``trivial'' lift, defined by setting $\lift(gr) = \theta(g)$, is always a valid choice of lift.

\begin{defn} 
Given a polarization $\theta\in \Gh$ and a good lift $\lift$ of $\theta$ (See Section~\ref{subsec:lifts}) and a  prestable LG-quasimap $\qmp = (\Ccal,x_1,\dots,x_k,\Pcal,\u,\spn)$, we say that $\qmp$ is \emph{$(\infty,\lift)$-stable} if 
\begin{enumerate}
\item There are no basepoints of $\u$ on $\Ccal$ and 
\item For every irreducible component $\Ccal'$ of $\Ccal$, the line bundle $\u^*(\Lcal_{\lift})$ restricted to $\Ccal'$ has nonnegative degree, with the degree only being allowed to vanish on components where $\klog$ is ample.
\end{enumerate}
\end{defn}

\begin{defn}
For a given choice of compatible $G$- and $\CC^*_R$-actions on a closed affine variety $Z\subseteq V$, a strongly regular character $\theta\in \Gh$ and a nondegenerate $W$,  we denote the corresponding stack of $k$-pointed, genus-$g$,
$0+$-stable LG-quasimaps into $\crst_\theta$ or $\X_\theta$ of degree $\beta$
by 
\[
\LGQ_{g,k}^{0+}(\crst_{\theta}, \beta) \quad \text{ or } \quad \LGQ_{g,k}^{0+}(\X_{\theta}, \beta),
\]
respectively.

If $\lift$ is a good lift of $\theta$, then the corresponding stack of $k$-pointed, genus-$g$,
$(\infty,\lift)$-stable LG-quasimaps into $\crst_\theta$ or $\X_\theta$ of degree $\beta$ is denoted by 
\[
\LGQ_{g,k}^{\infty,\lift}(\crst_{\theta}, \beta) \quad \text{ or }\quad \LGQ_{g,k}^{\infty,\lift}(\X_{\theta}, \beta),
\]
respectively.
\end{defn}

\section{A Hypersurface in Weighted Projective Space}

The remainder of this paper is dedicated to giving examples of the stack of stable LG-quasimaps to $\crst_\theta$ and $\X_\theta$ for various choices of input data.  We begin with an example of a {hypersurface} in weighted projective space.

\subsection{Basic Setup}

Suppose that $G = \CC^*$ and $F \in \CC[x_1,\dots,x_{K}]$ is a  quasihomogeneous polynomial of $G$-weights $(b_1, \dots, b_{{K}})$ and
 total $G$-degree $b$. Suppose further that $F$ has an isolated singularity at the origin.  Let
\[
W=p F{\colon} \CC^K \times \CC\rightarrow \CC,
\] 
where the variables $x_1,\dots, x_K$ are the first $K$ coordinates and $p$ is the last coordinate.
We assign $G$-weight $-b$ to the variable $p$, so that $W$ is $G$ invariant.  

The critical set of $W$ is given by the equations:
\[
\partial_p W=F=0 \dsand \partial_{x_i}W=p\partial_{x_i}F=0.
\]
If $p\neq 0$, then the fact that the only singularity of $F$ is at the origin means that $(x_1, \dots, x_{{K}})=(0, \dots, 0)$.  If 
$p=0$, then the only constraint is $F(x_1, \dots, x_{{K}})=0$.  So the critical locus is 
\[
\crit(W) = \{(0,\dots, 0, p) \mid p\in \C\} \cup \{(x_1,\dots, x_K,0) \in \C^K\times \C \mid  \ F(x_1,\dots, x_K) =0\}.
\]

Suppose that $b_i >0$ for $i=1, \dots, {{K}}$ and $b>0$.  If 
 $b = \sum_{i=1}^{{K}} b_i$, then we have a Calabi-Yau weight system, but we do not assume that here. 

\subsection{Phases for the Hypersurface}

Recall that different choices of linearization $\theta$ or moment map values $\mmapvalue$ give different quotients $[V\!\git{\theta} G]$, but any two linearizations lying in the same chamber have isomorphic quotients, called \emph{phases}.  

In the hypersurface case, the affine
 moment map
 $$\mmap=\frac12 \left(\sum_{i=1}^{{K}} b_i|x_i|^2- b |p|^2\right)$$ is a quadratic function
 whose only critical value is  $\mmapvalue=0$, and there are two phases $\mmapvalue>0$ or $\mmapvalue<0$.

\underline{Case of $\mmapvalue>0$:}

We have
 $$\sum_i^{{K}} b_i |x_i|^2=b|p|^2+2\mmapvalue.$$ For each choice of
 $p$, the set of $(x_1, \dots, x_{{K}}) \in \CC^{{K}}$, such that
 $(x_1,\dots, x_{K}, p) \in \mmap^{-1}(\mmapvalue)$, is a nontrivial
 ellipsoid $E$, isomorphic to $S^{2{{K}}-1}$; and we
 obtain a map from the symplectic quotient $\XS_\mmapvalue$ of $V$ to $ [E/U(1)] = \WP(b_1,\dots,b_{{K}})$.  
 
The resulting symplectic quotient $\XS_\mmapvalue$  can be expressed as the total space of the line  bundle $\Ocal(-b)$ over $\WP(b_1,\dots,b_{{K}})$.  If $\sum_i b_i = b$, this is the canonical bundle $\kk_{\WP(b_1, \dots,  b_{{K}})}$.

Alternatively, we can consider the GIT quotient $\left[\CC^{{{K}}+1}\git\theta
G\right]$, where $\theta:G\to \CC^*$ has weight $-\thetaweight$, 
with $\thetaweight>0$.  One can easily see
that the $\LV_\theta$-semistable points are $((\CC^{{K}}-\{\0 \})\times
\CC) \subset \CC^{{K}} \times \CC = \CC^{{{K}}+1}$, and the first projection
$pr_1:(\CC^{{K}}-\{\0 \})\times \CC \to (\CC^{{K}}-\{\0 \})$ induces the map
$[V\!\git{\theta} G] \to \WP(b_1,\dots,b_{{K}})$.

The critical locus $\crst_\theta = \{p=0 = F(x_1, \dots, x_{{K}})\}$ is a degree-$b$ hypersurface in the image of the zero section of $[V\!\git{\theta} G] \cong \Ocal(-d) \to \WP(b_1, \dots, b_{{K}})$. 
We call this phase the \emph{Calabi-Yau phase} or \emph{geometric phase}. In this setting, we write $\crst_{\mathrm{geom}}$ or $\X_{\mathrm{geom}}$ for $\crst_\theta$ or $\X_\theta$, respectively.

\underline{Case of $\mmapvalue<0$:}

We have
\[
\mmap^{-1}(\mmapvalue) = \left\{(x_1,\dots,x_{K},p) \,\middle |\, \sum_{i=1}^{{K}} b_i|x_i|^2-\mmapvalue= b |p|^2\right\}.
\]
For each choice of $x_1,\dots,x_{{K}} \in \CC^{{K}}$ the set of $p\in \CC$
such that $(x_1,\dots, x_{K}, p) \in \mmap^{-1}(\mmapvalue)$ is the circle
$S^1\subset \CC$, and
we obtain a map $\XS_\mmapvalue \to [S^1/U(1)]$.  If we choose the generator of
$U(1)$ to be $\lambda^{-1}$, then $p$ can be considered to have
positive weight $b$. Moreover, every $p$ has isotropy equal to
the $b$th roots of unity (isomorphic to $\ZZ_b$). The quotient
$[S^1/U(1)]$ is $\WP(b)=\B\ZZ_b = [\pt/\ZZ_b]$.
 
In the GIT formulation of this quotient, this corresponds to $\theta:G\to \CC^*$ of weight $-\thetaweight$,
and $\thetaweight<0$, the $\LV_\theta$-semistable points are equal to
$(\CC^{{{K}}} \times \CC^*) \subset \CC^{{{K}}+1}$. The second projection
$pr_2: (\CC^{{{K}}} \times \CC^*) \to \CC^*$ induces the map $[V\!\git{\theta} G]\to
\B\ZZ_b$.

The toric variety $\X_\theta = [V\!\git{\theta} G]$ can be viewed as the
total space of a rank-${{K}}$ orbifold vector bundle over $\B\ZZ_b$.  This bundle is actually just a $\ZZ_b$
bundle, where  $\ZZ_b$ acts by  
\[
(x_1, \dots, x_{{K}}) \mapsto (\xi_b^{b_1} x_1, \dots, \xi_b^{b_{{K}}} x_{{K}}) \qquad \xi_b = \exp(2\pi i/b).
\]
If we choose the $\CC^*_R$ action such that $W$ has $W$ has $\CC^*_R$-weight $b$, then this is exactly the
action of the element $J$ in FJRW-theory.  So the bundle $\X_\theta$ is
isomorphic to $[\CC^{{K}}/\genj ]$.  This is a special phase which is sort of like a toric variety with a finite group instead of
$\CC^*$.

 The critical locus is the single point $\{(0, \dots, 0)\}$ in the quotient $\X_\mmapvalue = [\CC^{{K}}/\ZZ_d]$.
It is clearly compact, so the the polynomial $W$ is nondegenerate.  We call $\X_\mmapvalue$ a \emph{Landau-Ginzburg phase} or a \emph{pure
  Landau-Ginzburg phase} \cite{Wit:97}. This example underlies
Witten's physical argument of the Landau-Ginzburg/Calabi-Yau
correspondence for Calabi-Yau hypersurfaces of weighted projective
spaces.

In this setting, we write $\crst_{\mathrm{LG}}$ or $\X_{\mathrm{LG}}$ for $\crst_\theta$ or $\X_\theta$, respectively.

\subsection{GLSM Moduli Space for the Hypersurface}

\subsubsection{Geometric phase}

For the phase $\mmapvalue>0$, we choose  $\CC^*_R$-weights $\Rq_{x_i}=0$ and $\Rq_{p}=1$ (thus giving  a {hybrid model}), so $W$ has $\CC^*_R$-weight $d=1$.  The element $J$ is trivial, and the group 
\[
\Gamma = \{(g^{b_1}, \dots, g^{b_K}, g^{-b}r) \mid g\in G, r\in \CC^*_R\} 
\]
 is a direct product $\Gamma \cong G \times \CC^*_R$, with $\spl$ and $\chiR$  just the first and second projections, respectively. 

There are two ingredients to the moduli space: an Artin stack of geometric data and a more subtle stability condition.  

The geometric data correspond to sections of the vector bundle $\Pcal\times_\Gamma V$.  This bundle can be written as a direct sum of line bundles 
\[
\Ecal = \Pcal\times_\Gamma V \cong \Acal^{\otimes b_1} \oplus \Acal^{\otimes b_2} \oplus \cdots \oplus \Acal^{\otimes b_K} \oplus (\Acal^{\otimes(-b)}\otimes\Bcal),
\]
where $\Acal$ corresponds to the $G$-action (the map $\Ccal \rTo^{\Pcal} B\Gamma \rTo^{\spl} BG$)
and $\Bcal$ corresponds to the $\C^*_R$ action (the map $\Ccal \rTo^{\Pcal} B\Gamma \rTo^{\chiR} \CC^*_R$). And we have an isomorphism $\spn: \chiR_*\Pcal = \Bcal \to \klogc$.  Therefore, the section $\u$ corresponds to a sequence $\u = (s_1,\dots, s_K, p)$ where $s_i$ is a section of $\Acal^{b_i}$ and $p$ is a section of $\Acal^{-b}\otimes \klogc$.

So LG-quasimaps to the geometric phase $\X_{\mathrm{geom}}$  correspond to the following data:
\begin{equation}\label{eq:toric-geom-phase-X}
      \{(\Ccal, \Acal, s_1, \cdots, s_{K},p) \mid  s_i\in H^0(\Ccal,\Acal^{b_i}),\, p\in H^0(\Ccal,\Acal^{-b}\otimes \klogc) \}.
      \end{equation}

For all stability conditions, the LG-quasimaps can fail in at most a finite number of points to map to the $\theta$-semistable locus.  This corresponds to the locus where $s_1,\dots,s_K$ do not simultaneously vanish.  Moreover, $\crit(W)^{ss}_{\mathrm{geom}}$ has $p=0$.  But if $p=0$ at all but a finite number of points, then $p$ is always $0$.

Thus, without further specifying stability conditions, LG-quasimaps to 
the geometric phase 
$\crst_{\mathrm{geom}}$ correspond to the following data
\begin{equation}\label{eq:toric-geom-phase-CR}
      \{(\Ccal, \Acal, s_1, \cdots, s_{K}) \mid  s_i\in H^0(\Ccal,\Acal^{b_i})\},
\end{equation}
where $\Ccal$ is a marked orbicurve, $\Acal$ is a line bundle over $\Ccal$, and the section $\u = (s_1,\dots,s_K,0)$ maps to $\Pcal\times_\Gamma \crit(W)$.

Now we consider the stability conditions.

\underline{Case of $\ve = \infty$:} 

We must find a good lift of $\theta$.   Let $\ell$ be a generator of $\LV_\theta^*$ over $\CC[V^*]$ with $G$ acting on $\ell$ with weight $-\thetaweight$ (and $\thetaweight >0$ in the geometric phase).   The trivial lift $\lift_0$ of $\theta$  corresponds to $\CC^*_R$ acting trivially on $\ell$.  A monomial of the form  $x_i^{\thetaweight}\ell^{b_i}$ is $\Gamma$-invariant and does not vanish on points with $x_i\neq 0$, so every point of $\CC^{{N}}\times
\CC$ with $x_i\neq 0$ is in $V^{ss}_\Gamma(\lift_0)$.   Letting $i$ range from $1$ to ${N}$ shows that $V^{ss}_\Gamma(\lift_0) = V^{ss}_G(\theta)$. Thus  $\lift_0$ is a good lift of the character $\theta$.  

Any other lift $\lift$ must have nontrivial $\CC^*_R$ action on $\ell$ and thus any $\C^*_R$-invariant element of $\C[x_1,\dots, x_K,p][\ell]$ must have each monomial containing a power of $p$, and hence points with $p=0$ will not be $\lift$-semistable.  Therefore, $\lift_0$ is the only good lift of $\theta$.

Finally, $\sigma^*\Lcal_{\lift_0}$ is determined by the action of $\Gamma$ on $\LV_\theta$ (or the inverse of the action of $\lift_0$ on $\ell$), so in this case $\sigma^*\Lcal_{\lift_0}  = \Acal^{\thetaweight}$. Thus, the $(\infty,\lift_0)$-stable LG-quasimaps to $\crst_{\mathrm{geom}}$ or to $\X_{\mathrm{geom}}$ consist of those data \eqref{eq:toric-geom-phase-CR} or \eqref{eq:toric-geom-phase-X}, respectively, 
satisfying the conditions that
\begin{enumerate}
\item The section $\u$ has no basepoints (the $s_i$ never vanish simultaneously).
\item The line bundle $\Acal$ has positive degree on every component $\Ccal'$ of $\Ccal$ where $\klogc$ has nonpositive degree.
\end{enumerate}

Thus $(\infty,\lift_0)$-stable LG-quasimaps to $\crst_{\mathrm{geom}}$ correspond  to stable maps to the hypersurface $X_F = \{F=0\} \subset W\PP(b_1,\dots, b_K)\}$. And $(\infty,\lift_0)$-stable LG-quasimaps to $\X_{\mathrm{geom}}$ are stable maps to $X_F$ with $p$-fields, studied in \cite{ChaLi:11, CLL:13}.

\underline{Case of $\ve = 0+$:} 

The $0+$ stable LG-quasimaps must have the section $\u$ take its values in $\crit(W)^{ss}_{\mathrm{geom}}$ for all but a finite number of basepoints $\mrkp\in \Ccal$.  

Letting $\lift$ be the trivial lift of $\theta$, we have $\sigma^*\Lcal_\lift = \Acal^{\thetaweight}$ and so the $0+$ stable LG quasimaps to $\crst_{\mathrm{geom}}$ or to $\X_{\mathrm{geom}}$ 
are the data of \eqref{eq:toric-geom-phase-CR} or \eqref{eq:toric-geom-phase-X}, respectively, satisfying the stability conditions:
\begin{enumerate}
\item There are at most a finite number of basepoints (where the $s_i$ vanish simultaneously), and these only occur away from nodes and marked points.
\item Every rational component has at least two special points (node or marked point).
\item On rational components with exactly two special points, the line bundle $\Acal$ has positive degree.
\end{enumerate}

Thus $0+$-stable LG-quasimaps into $\crst_{\mathrm{geom}}$ are stable quotients into $X_F \subset \WP(b_1,\dots,b_{K})$, and $0+$-stable LG-quasimaps into $\X_{\mathrm{geom}}$ are stable quotients into $X_F$ with $p$-fields.

\begin{rem}
There is a parallel theory of quasimaps  into $X_F$. Both theories have the same moduli spaces, but the virtual cycle constructions are different. For $\ve=\infty$, Chang-Li \cite{ChaLi:11} 
      proved the equivalence using a sophisticated degeneration argument. A similar argument probably works for other $\ve$-theories.
\end{rem}

\subsubsection{LG phase}

For the phase $\mmapvalue<0$, we choose $\CC^*_R$ to have weights
 $\Rq_{x_i}=b_i$ and $\Rq_{p}=0$, which again gives a hybrid model. 
But now $W$ has $\CC^*_R$-weight $d=b$, and 
 $J = (\xi^{b_1}, \dots, \xi^{b_{K}},1)$, where $\xi = \exp(2\pi i /d)$. We have 
 \begin{align}
 \Gamma &= \{((gr)^{b_1}, \dots, (gr)^{b_{K}}, g^{-b}) \mid g\in G, r\in \CC^*_R\} \label{eq:toric-LG-Gamma1}  \\ 
 &= \{(\alpha^{b_1}, \dots, \alpha^{b_{K}}, \beta) \mid \alpha,\beta \in \CC^*\},\label{eq:toric-LG-Gamma2}
 \end{align}
where $\alpha = gr$ and $\beta = g^{-b}$. 
  with $\chiR:\Gamma \to \CC^*$ given by $(\alpha^{b_1}, \dots, \alpha^{b_{K}}, \beta) \mapsto \alpha^b\beta$.  

Thus, the vector bundle $\Ecal = \Pcal\times_\Gamma V$ associated to an LG-quasimap is a direct sum of line bundles on $\Ccal$:
\[
\Ecal = \Acal^{b_1}\oplus\cdots\oplus\Acal^{b_{K}}\oplus\Bcal, 
\] 
where $\Acal$ corresponds to $\alpha$ and $\Bcal$ corresponds to $\beta$ in the presentation \eqref{eq:toric-LG-Gamma2} of $\Gamma$.  Moreover, we have $\spn: \Acal^b\otimes\Bcal \to \klogc$ is an isomorphism.  

Thus, LG-quasimaps to $\X_{\mathrm{LG}}$ again consist of exactly the same data as \eqref{eq:toric-geom-phase-X}:
   \[
      \{(\Ccal, \Acal, s_1, \cdots, s_{K}, p)\mid s_i\in H^0(\Ccal,\Acal^{b_i}),\  p\in H^0(\Acal^{-b}\otimes \klogc)\}.
      \]
The base points of these quasimaps occur precisely at the zeros of $p$, and the base locus forms an effective divisor $D$ in $\Ccal$ with  $\Acal^{-b}\otimes \klogc\cong  O(D)$, so the section $p$ gives an isomorphism $\Acal^b \cong \klogc(-D)$ and can be viewed as a \emph{weighted $b$-spin} condition (see \cite{RR}).
So we can reformulate the moduli problems as 
  \begin{equation}\label{eq:toric-LG-phase-X}
      \{(\Ccal, \Acal, D, s_1, \cdots, s_{K}) \mid s_i\in H^0(\Ccal,\Acal^{b_i}),  \Acal^b \cong \klogc(-D)\},
      \end{equation}
where each $D$ is an effective divisor that is disjoint from the nodes and marked points of $\Ccal$.

For the LG-quasimaps to lie in the critical locus requires every $s_i=0$, so LG-quasimaps to $\crst_{\mathrm{LG}}$ consist of
  \begin{equation}\label{eq:toric-LG-phase-CR}   
      \{(\Ccal, \Acal, D) \mid  \Acal^b \cong \klogc(-D)\}.
    \end{equation}
    
\underline{Case of $\ve = \infty$:} 

Again, the trivial lift is the only good lift of $\theta$. To see this, let $\ell$ be a generator of $\LV_\theta^*$ over $\CC[V^*]$ with $G$ acting on $\ell$ with weight $-\thetaweight$ (and $\thetaweight<0$ in the LG phase). The trivial lift $\lift_0$ corresponds to $\CC^*_R$ acting trivially on $\ell$, and a monomial of the form  $p^{-\thetaweight}\ell^b$ is $\Gamma$-invariant and does not vanish on points with $p\neq 0$, so $V^{ss}_\Gamma(\lift_0) = V^{ss}_G(\theta)$. Thus  $\lift_0$ is a good lift of the character $\theta$.  Any other lift $\lift$ must have nontrivial $\C^*_R$-action on $\ell$, and hence any $\Gamma$-invariant function must have at least one factor of $x_i$ in every monomial, which implies that any point of $V$ with $x_1=x_2=\cdots=x_k =0$ is not $\lift$-semistable.  Thus $\lift_0$ is the only good lift.  

The line bundle $\sigma^*\Lcal_{\lift_0}$ is determined by the action of $\Gamma$ on $\LV_\theta$, that is by $g^e$ in the presentation \eqref{eq:toric-LG-Gamma1}, which implies that 
\[
\sigma^*\Lcal_{\lift_0} \cong \Bcal^{-e/b} \cong \klogc^{-e/b} \otimes \Acal^{e}.
\]
Since $\ve =\infty$, no base points are permitted, so $D=0$ and $\Acal^b \cong \klogc$.  For convenience, let us assume that $e = -cb$ for some $c>0$.  The stability condition is now that 
\[
\klogc^{c} \otimes \Acal^{-cb} \cong \klogc^{c} \otimes \klogc^{-c} = \Ocal
\]
can only have degree $0$ on components where $\klogc$ is ample, thus $\Ccal$ must be a stable orbicurve.

So in this case $(\infty,\lift_0)$-stable LG-quasimaps to $\crst_{\mathrm{LG}}$ correspond to stable $b$-spin curves
\[
\{ \Ccal, \Acal \mid \Acal^b \cong \klogc\},
\]
studied in \cite{JKV, AbJa:03}.

\underline{Case of $\ve = 0+$:} 

In this case basepoints are permitted, so $D$ is not necessarily $0$. 
The bundle
\[
\sigma^*\Lcal_{\lift_0} \cong \klogc^{c} \otimes \Acal^{-bc} \cong \Ocal(cD)
\] 
must have positive degree on any component where $\klogc$ is not ample.  

Thus $0+$-stable LG quasimaps to $\crst_{\mathrm{LG}}$ correspond to the data of \eqref{eq:toric-LG-phase-CR} satisfying the conditions:
\begin{enumerate}
\item Every rational component has at least two special points, and 
\item On every irreducible component $\Ccal'$ with trivial $\omega_{log,\Ccal'}$, there is at least one basepoint.
\end{enumerate}

\section{Complete Intersection in Weighted Projective Space}

\subsection{Basic Setup}

Suppose that $G = \CC^*$ and we have several quasihomogeneous polynomials $F_1, F_2,
\dots, F_M \in \CC[x_1,\dots,x_{{K}}]$ of $G$-degree $(d_1, \dots, d_M)$, where each variable $x_i$ has $G$-weight $b_i>0$. We assume that the $F_j$
 intersect transversely in $\WP(b_1, \dots, b_K)$ and define a
 complete intersection. Let
 \[
W=\sum_i p_i F_i{\colon} \CC^{{{K}}+M}\rightarrow \CC,
\]
 where we assign
 $G$-weight $-d_i$ to $p_i$. In the special case that 
 $\sum_i b_i=\sum_j d_j$,
 then the complete intersection defined by $F_1=\cdots=F_M=0$
 is a Calabi-Yau orbifold in $\WP(b_1,\dots,b_{{K}})$.  But we do not assume the  Calabi-Yau condition here. 
 
The critical set of $W$ is defined by the following equations:
\begin{equation}\label{eq:crit-for-toricLG}
\partial_{p_j}W=F_j=0, \ \partial_{x_i}W=\sum_j p_j \partial_{x_i}
 F_j=0.
\end{equation} 
Since the $F_j$ intersect transversely, an argument similar to that for the hypersurface shows that the critical locus consists of those $(\x,\p)$ where either $\p=\0$ and $\x$ satisfies $F_i(\x)=0$ for all $i$, or $\p$ is unconstrained and $\x=\0$:
\[
\crit(W) = \{(\0,\p)\in \C^K\times \C^M \mid \p \in \C^M \}\cup \{ (\x,\0) \in \C^K\times \C^M \mid F_i(\x) = 0\, \forall i\}.
\]

\subsection{Phases for a Complete Intersection}
 The moment map is 
\[
  \mmap=\sum_i  \frac12 b_i|x_i|^2-\frac12 \sum_j d_j|p_j|^2.
\] Again, there are two phases,
$\mmapvalue>0$ and $\mmapvalue<0$. 

\underline{Case of $\mmapvalue>0$:}

When $\mmapvalue>0$, we again call this the \emph{geometric phase}.
Any choice of $\p = (p_1,\dots
p_M)$ determines a nontrivial ellipsoid $E \subset \CC^{{K}}$ of 
points $\x = (x_1, \dots, x_{{K}})$ such that
  $(\x,\p)$ lies in $\mmap^{-1}(\mmapvalue)$.  
Quotienting by $U(1)$, the first projection $pr_1:E \times \CC^M \to
  E$ induces a map $\X_{\mathrm{geom}} \to \WP(b_1, \dots, b_{{K}})$.
The full quotient is $\X_{\mathrm{geom}} = \bigoplus_j \Ocal(-d_j)$ over
 $\WP(b_1, \dots, b_{{K}})$.

In the GIT formulation, this again corresponds to $\theta:G\to \CC^*$ having weight $-\thetaweight$, with $\thetaweight>0$.  The semistable points of this phase are those with $\x \neq \0$,
and the semistable points of the critical locus correspond to points in 
 \[
 \{F_1 = \dots = F_M = 0\} 
 \]
The quotient $\crst_{\mathrm{geom}}$ is the locus of the complete intersection defined by all the $F_i$ vanishing in the zero section of $\X_{\mathrm{geom}} \to \WP(b_1,\dots,b_{{K}})$.

\underline{Case of $\mmapvalue<0$:}

When $\mmapvalue < 0$, we again call this the \emph{Landau-Ginzburg phase}.  The quotient is $\X_{\mathrm{LG}} = \bigoplus_i
  \Ocal(-b_i)$ over $\WP(d_1, \dots, d_M)$.
 
In the GIT formulation, this corresponds to $\theta:G\to \CC^*$ having weight $-\thetaweight$, with $\thetaweight<0$. The semistable points of this phase are those with $\p \neq \0$,
and the semistable points of the critical locus correspond to the image of the zero section of $\X_{\mathrm{LG}} \to \WP(d_1,\dots,d_M)$.

\subsection{GLSM Moduli Space for a Complete Intersection}

\subsubsection{Geometric Phase}

We choose the R-charge to act on $\C^K\times \C^M$ with weights $(0,\cdots, 0, 1,\cdots, 1)$, which gives a hybrid model.  And $W$ has $\CC^*_R$-weight $d=1$.  The element $J$ is trivial, and the group 
\[
\Gamma = \{(g^{b_1}, \dots, g^{b_K}, g^{-d_1}r, \dots g^{-d_M}r) \mid g\in G, r\in \CC^*_R\} 
\]
 is again a direct product $\Gamma \cong G \times \CC^*_R$, with $\spl$ and $\chiR$  just the first and second projections, respectively. 

The geometric data correspond to sections of the vector bundle $\Ecal = \Pcal\times_\Gamma V$, which can be written as a direct sum of line bundles 
\[
\Ecal =  \Acal^{b_1} \oplus \Acal^{b_2} \oplus \cdots \oplus \Acal^{b_K} \oplus (\Acal^{-d_1}\otimes\Bcal) \oplus \cdots \oplus (\Acal^{-d_M}\otimes\Bcal),
\]
where $\Acal$ corresponds to the $G$-action and $\Bcal$ corresponds to the $\C^*_R$ action.  And we have an isomorphism $\spn: \chiR_*\Pcal = \Bcal \to \klogc$.

So LG-quasimaps to the geometric phase $\X_{\mathrm{geom}}$  correspond to the data:
\begin{equation}\label{eq:CI-geom-phase-X}
      \{(\Ccal, \Acal, s_1, \cdots, s_{K},p_1,\dots, p_M) \mid  s_i\in H^0(\Ccal,\Acal^{b_i}),\, p_i\in H^0(\Ccal,\Acal^{-d_i}\otimes \klogc) \}.
      \end{equation}

Again $\crit(W)^{ss}_{\mathrm{geom}}$ has $\p=0$, so without further specifying stability conditions, LG-quasimaps to 
the geometric phase 
$\crst_{\mathrm{geom}}$ correspond to the following data
\begin{equation}\label{eq:CI-geom-phase-CR}
      \{(\Ccal, \Acal, s_1, \cdots, s_{K}) \mid  s_i\in H^0(\Ccal,\Acal^{b_i})\},
\end{equation}
where $\Ccal$ is a marked orbicurve, $\Acal$ is a line bundle over $\Ccal$, and the section $\u = (s_1,\dots,s_K,\0)$ maps to $\Pcal\times_\Gamma \crit(W)$.

\underline{Case of $\ve = \infty$:} 

Again $\lift_0$  is easily seen to be a good lift, and $\sigma^*\Lcal_{\lift_0}  = \Acal^{\thetaweight}$. Thus, just as in the hypersurface case, we have that $(\infty,\lift_0)$-stable LG-quasimaps to $\crst_{\mathrm{geom}}$ correspond  to stable maps to the complete intersection  $X = \{F_1=\cdots = F_M = 0\} \subset W\PP(b_1,\dots, b_K)\}$. And $(\infty,\lift_0)$-stable LG-quasimaps to $\X_{\mathrm{geom}}$ are stable maps to $X$ with $p$-fields.

\underline{Case of $\ve = 0+$:} 

The arguments given in the hypersurface case are easily adapted to show that $0+$-stable LG-quasimaps into $\crst_{\mathrm{geom}}$ are stable quotients into $X \subset \WP(b_1,\dots,b_{K})$, and $0+$-stable LG-quasimaps into $\X_{\mathrm{geom}}$ are stable quotients into $X_F$ with $p$-fields.

\subsubsection{LG phase}

Assume that $d_1=\cdots=d_r=d$ and choose the R-charge weights $\Rq_{x_i}=b_i$ and $\Rq_{p_j}=0$.  
Now $W$ has $\CC^*_R$-weight $d$, and 
 \begin{align}
 \Gamma &= \{ (gr)^{b_1}, \dots, (gr)^{b_K}, g^{-d}, \dots, g^{-d} \mid g\in G, r\in \C^*_R\}  \label{eq:CI-LG-Gamma1}  \\ 
 & =\{ \alpha^{b_1}, \dots, \alpha^{b_K}, \beta, \dots, \beta) \mid \alpha,\beta \in \C^*\}, 
\label{eq:CI-LG-Gamma2}
 \end{align}
 where $\alpha = gr$ and $\beta = g^{-d}$, and the map $\chiR$ sends $(\alpha,\beta)$ to $\alpha^d\beta$.
 
Thus, $\Ecal = \Pcal\times_\Gamma V$ is a direct sum of line bundles on $\Ccal$:
\[
\Ecal = \Acal^{b_1}\oplus\cdots\oplus \Acal^{b_{K}}\oplus\Bcal\oplus \cdots \oplus \Bcal, 
\] 
where $\Acal$ corresponds to $\alpha$ and $\Bcal$ corresponds 
to $\beta$ in the presentation \eqref{eq:CI-LG-Gamma2} of $\Gamma$.  Moreover, $\spn: \Acal^d\otimes\Bcal \to \klogc$ is an isomorphism.  

Thus, LG-quasimaps to $\X_{\mathrm{LG}}$ consist of the data:
   \[
      \{(\Ccal, \Acal, s_1, \cdots, s_{K}, p_1,\dots, p_M)\mid s_i\in H^0(\Ccal,\Acal^{b_i}),\  p_i\in H^0(\Acal^{-d}\otimes \klogc)\}.
   \]
And LQ-quasimaps to $\crst_{\mathrm{LG}}$ also require that all the $s_1,\dots, s_K$ vanish, giving
\[
      \{(\Ccal, \Acal, p_1,\dots, p_M)\mid  p_i\in H^0(\Acal^{-d}\otimes \klogc)\}.
 \]

\underline{Case of $\ve = \infty$:} 

Again, the trivial lift is the only good lift of $\theta$. The line bundle $\sigma^*\Lcal_{\lift_0}$ is again
\[
\sigma^*\Lcal_{\lift_0} \cong \Bcal^{-e/d} \cong \klogc^{-e/d} \otimes \Acal^{e}.
\]
And so the stability condition is that 
\[
\klogc \otimes \Acal^{-d} 
\]
can only have degree $0$ on components where $\klogc$ is ample.
And since $\ve =\infty$, no base points are permitted, so the $p_i$ cannot all simultaneously vanish.

So in this case $(\infty,\lift_0)$-stable LG-quasimaps to $\crst_{\mathrm{LG}}$ correspond to stable maps to $W\PP(d,d,\dots,d)$. 
And for each  
$(p_1, \dots, p_n) \in \WP(d,\dots, d)$, we have a pure 
LG-model of superpotential
$\sum_i p_i F_i$.  One can view this as a family of pure LG-theories. 

\begin{rem}\label{rem:infty-all-same}
An LG-phase of a complete intersection of differing degrees (i.e., $d_i \neq d_j$ for some $1\le i,j \le r$) does not admit a hybrid model structure and will generally have no good lift.
\end{rem}

%% BOOKMARK

\underline{Case of $\ve = 0+$:} 

Now basepoints are permitted, and every rational component must have at least two special points.  Again, the stability condition is that 
\[
\klogc \otimes \Acal^{-d} 
\]
can only have degree $0$ on components where $\klogc$ is ample, 

\begin{rem}
As mentioned in Remark~\ref{rem:infty-all-same}, without the condition $d_1=d_2 = \cdots = d_r=d$ there is usually no good lift of $\theta$.  But in the $\ve=0+$ case, we do not need a good lift, so fixing any $d>0$ we can take $\Rq_{p_j}=d-d_j$, which again gives $W$ the  $\CC^*_R$-weight of $d$ and a corresponding proper DM stack of $0+$-stable LG-quasimaps. 
\end{rem}

\section{Graph Moduli Space}

  The graph moduli space is very important in Gromov-Witten theory.   We can construct it in the GLSM setting as follows.

Suppose that we have a phase $\theta$, a superpotential $W: [\CC^n/G]\to  \CC$ with a certain $R$-charge $\CC^*_R$, defining $\Gamma$ and a lift $\lift$ of $\theta$. We construct a new GLSM as follows.  

Let $V' = V\times \CC^2$, and let $\CC^*$ act on $\CC^2$ with weights $(1,1)$. Let $G' = G\times \CC^*$ act on $V'$ with the product action, so $G$ acts trivially on the last two coordinates and $\CC^*$ acts trivially on the first $n$ coordinates.

Let $\theta':G'\to \CC^*$ be given by sending any $(g,h) \in G\times \CC^*$ to $\theta(g) h^{-k}$ for some $k>0$.  The GIT quotient is the product $[V'\git{\theta'}G'] = [V\git{\theta}G]\times \PP^1$.  Let $W'$ be defined on $V'$ by the same polynomial as $W$, so that the critical locus of $W'$ in $V'$ is $\CC^2$ times the the critical locus of $W$, and the GIT quotient of the critical locus is the product of $\PP^1$ and the corresponding quotient in the original GLSM.

Keeping the same $R$-charge (that is, letting $\CC^*_R$ acts trivially on the last two coordinates of $V'$), we have $\Gamma' = \Gamma \times \CC^*$, and we construct a  lift $\lift'$ of $\theta'$ by sending $(\gamma,h) \in \Gamma \times \CC^*$ to $\lift(\gamma)h^{-k}$. It is easy to see that $\lift'$ is a good lift of $\theta'$ if $\lift$ is a good lift of $\theta$.

In the $\ve = \infty$ case, no basepoints can occur, and projecting to the two new coordinates $(z_0,z_1)$  induces a stable map $\Ccal\rightarrow \PP^1$. Therefore, the new GLSM in this case can be reformulated as the usual 
GLSM for $[V\git{\theta}G]$ with the additional 
data of a stable map $f: \Ccal \rightarrow \PP^1$.

  \section{Generalized  graph space}\label{exa:RRS}

  We can generalize slightly the graph moduli space to obtain a new moduli space with a remarkable property. Let's take the quintic GLSM as an example. Now, we consider a new GLSM on $\CC^{6+2}\git{}(\CC^*)^2$, given by $G = (\CC^*)^2$ acting on $V = \CC^8$ with weights
  $$\left(\begin{array}{cccccccc}
  1&1&1&1&1&-5&d&0\\
  0&0&0&0&0&0&1&1
  \end{array}\right)$$
  for an integer $d>0$. 
Let the coordinates on $V$  be $x_0,\dots, x_4, p, z_0, z_1$, corresponding the the columns in the weight matrix above.  The moment maps are
  $$\mu_1=\frac{1}{2}\left(\sum_{i=0}^4 |x_i|^2-5|p|^2+d|z_0|^2\right),\     \mu_2=\frac{1}{2}(|z_0|^2+|z_1|^2).$$
  There are three chambers. We are interested primarily in the chamber $0<\mu_1<d\mu_2$.  This corresponds to a character $\theta$ of $G$ with weights $(-e_1,-e_2)$ and $0<e_1<d e_2$.  The $\theta$-unstable locus for $\theta$ is 
  \[
  \{x_0=x_1=x_2=x_3=x_4=z_0= 0\}\cup \{p=z_1=0\}\cup\{z_0=z_1= 0\}.
  \]

Taking the superpotential $W = \sum_{i=1}^5 p x_i^5$ and the $R$-charge of weight $(0,0,0,0,\linebreak[1] 0,1,0,0)$, we have $\Gamma = G\times \CC^*_R = \{(a,a,a,a,a,ra^{-5},ba^d,b)\mid 
a,b,r\in \CC^*\}$, and the map $\chiR$ takes  $(a,a,a,a,a,ra^{-5},ba^d,b)$ to $r$.

There is no good lift of $\theta$, so we restrict to the case of $\ve = 0+$.  We must choose some lift for the stability condition, so we take the trivial lift $\lift(a,b,r) = a^{-e_1}b^{-e_2}$.  Any other lift will give the same stability conditions.

The resulting moduli problem consists of   
\begin{align*}
 \{(\Ccal, \Acal, \Bcal, x_0, \cdots, x_4, p, z_1,z_2)\mid   x_i\in &H^0(\Ccal,\Acal), p\in H^0(\Ccal,\Acal^{-5}\otimes \klogc) \\
 & z_1\in H^0( \Ccal,\Acal^{d}\otimes \Bcal), 
  z_2\in H^0( \Ccal,\Bcal)\},
  \end{align*}
satisfying the stability condition that $\u^*\Lcal_{\lift} =\Acal^{-e_1}\Bcal^{-e_2}$ is ample on all components where $\klogc$ has degree $0$.

The critical locus is 
\[
\crit(W)=\{x_0=\cdots = x_4 = 0\} \cup \{ p = \sum_{i=0}^4 x_i^5=0\}
\]
And so, in the chamber we are interested in, the GIT quotient of the critical locus has two components $\crst_{\theta} = C_1 \cup C_2$, where 
$C_1$ corresponds to $x_0=\dots = x_4 =0$, so $z_0\neq 0$ can be scaled to $1$ by the second $\C^*$ action and so the quotient $C_1$ is isomorphic to $\PP(5,1)$, with coordinates $p,z_1$.  The component $C_2$ corresponds to the locus $\{ p = \sum_{i=0}^4 x_i^5=0\}$

The critical locus admits a $\CC^*$ action by multiplication on $z_1$.  The fixed loci of the action are 
\begin{enumerate}
\item  The locus $\{x_0=\dots=x_4=z_1=0\}$, which is the point $B\mu_5$ inside $\PP(5,d)$.  
\item The locus $\{p = z_0 = \sum_{i=0}^4 x_i^5 = 0\}$, which is a Calabi-Yau threefold in $\PP^4$.  
\item  The locus $\{x_0=\dots=x_4=p=0\}$, which is the point $B\mu_d$ inside $\PP(5,d)$.
\end{enumerate}

The GLSMs for the three fixed loci correspond, respectively to (1) a weighted FJRW theory,  (2) $0+$-stable quasimaps to the quintic threefold, and (3)  
the theory of Hassett stable curves
with light points given by the vanishing of $z_0$ and $z_1$.
This remarkable property gives us the hope that we can extract a relation between Gromov-Witten theory and FJRW-theory
 geometrically by using localization techniques on this moduli space. A program is being carried out right now for the $\ve=0^+$ theory  \cite{RRS,CJR}.

%If $z_1=0$ and $p$ never vanishes.  In this case $p$ gives an isomorphism $\Acal^5\cong \klogc$, and the moduli problem is   
%  \[
% \{(\Ccal, \Acal,\Bcal, z_0) \mid \Acal^5\cong \klogc, z_0 \in H^0(\Ccal,\Acal^{d}\otimes \Bcal) \},
%\]
%    with the stability condition that $\Bcal^{-e_2}$ is ample whenever $\klogc$ has degree $0$.
%  
%  
%  
%  If z_1 never vanishes, it gives an isomorphism $\Bcal \cong \Acal^{-d}(-D)$ for some effective divisor $D$ and the resulting moduli problem is 
%\[
% \{(\Ccal, \Acal, x_1, \cdots, x_5, p)\mid   x_i\in &H^0(\Ccal,\Acal), p\in H^0(\Ccal,\Acal^{-5}\otimes \klogc) \},
%\]
%satisfying the stability condition that $\Acal^{-e_1+de_2}(-e_2 D)$ is ample on all components where $\klogc$ has degree $0$.  
%

  The same theory with a different $\ve=\infty$ stability condition was discovered and the localization argument was carried out independently by Chang-Li-Li-Liu \cite{CLL:15}.

\section{Non-Abelian examples}
The subject of gauged linear sigma models for non-Abelian groups is a
very active area of research in physics and is far from complete.
Here, we discuss complete intersections in a Grassmannian or flag variety.

All of this should work similarly in the
setting of complete intersections of quiver varieties, although the
details have not been worked out. It would be very interesting to explore  
mirror symmetry among Calabi-Yau complete intersections in quiver
varieties.

\subsection{Complete Intersections in a Grassmanian} 

The space $\Gr(k, {n})$ can be constructed as a GIT  quotient $M_{k,{n}}\git{}\GL(k, \CC)$, where $M_{k,{n}}$ is the space of
$k\times {n}$ matrices and $\GL(k, \CC)$ acts as matrix multiplication on
the left.  

The Grassmannian $\Gr(k, {n})$ can also be embedded into
$\PP^{{{K}}}$ for ${{K}}=\frac{{n}!}{k!({n}-k)!}-1$ by the Pl\"ucker embedding
$$A \mapsto (\dots, \det(A_{i_1,\dots, i_k}), \dots),$$ where
$A_{i_1, \cdots, i_{k}}$ is the $(k\times k)$-submatrix of $A$ consisting of the columns $i_1, \dots,
i_k$. 

The group $G = \GL(k,\CC)$ acts on the Pl\"ucker coordinates $B_{i_1, \cdots, i_k}(A)=\det(A_{i_1,\cdots, i_k})$ by the determinant, that is, for any $U\in G$, and $A\in M_{k,{n}}$ we have
$$B_{i_1, \cdots, i_k}(UA)=\det(U)B_{i_1, \cdots, i_k}(A),$$

Let $F_1, \dots, F_s\in \CC[B_{1,\dots,k}, \dots, B_{{n}-k+1,\dots,{n}}]$ be
degree-$d_j$ homogeneous polynomials such that the zero loci $Z_{F_j}=\{F_j=0\}$ and the Pl\"ucker embedding of  $\Gr(k,{n})$ all intersect transversely in $\PP^{{{K}}}$.  We let 
$$Z_{d_1, \cdots, d_s}=\Gr(k,{n})\cap \bigcap_{j} Z_{F_j}$$
denote the corresponding complete intersection.
 
The analysis of $Z_{d_1, \cdots, d_s}$ is similar to the Abelian case. Namely, let
$$W=\sum_j p_j F_j{\colon} M_{k,{n}}\times \CC^s\rightarrow \CC$$
be the superpotential.
We assign an action of $G = \GL(k, \CC)$ on $p_j$ by $p_j\rightarrow
\det(U)^{-d_j}$.

The phase structure is similar to that of a complete intersection in projective space.  The moment map is given by $\mmap(A,p_1,\dots,p_s) = \frac{1}{2} (A\bar{A}^T - \sum_{i=1}^s d_i |p_i|^2)$.  Alternatively, to construct a linearization for GIT, the only characters of $GL(k,\CC)$ are powers of the determinant, so $\theta(U) = \det(U)^{-\thetaweight}$ for some $\thetaweight$, and $\mmapvalue$ will be positive precisely when $\thetaweight$ is positive.

Let $\ell$ be a generator of $\CC[\LV^*_\theta]$ over $\CC[V^*]$.   Any element of $H^0(V,\LV_\theta)$ can be written as a sum of monomials in the Pl\"ucker coordinates $B_{i_1,\dots,i_k}$ and the $p_j$ times $\ell$.  Any $U\in G$ will act on a monomial of the form $\prod B_{i_1,\dots,i_k}^{b_{i_1,\dots,i_k}} \prod p_j^{a_j}\ell^m$ by multiplication by $\det(U)^{\sum b_{i_1,\dots,i_k} -\sum d_j a_j -m\thetaweight}$.  

\subsubsection{Geometric Phase}

Assume that $\thetaweight>0$.
In order to be $G$-invariant, a monomial must have $\sum b_{i_1,\dots,i_k}>0$, which implies that any points with every $B_{i_1,\dots,i_k} =0$ must be unstable, but for each $m>0$ and each $k$-tuple $(i_1,\dots,i_k)$ the monomial $B_{i_1,\dots,i_k}^{m\thetaweight} \ell^m$ is $G$ invariant, so every point with at least one nonzero $B_{i_1,\dots,i_k}$ must be $\theta$-semistable.  
Thus $[V\!\git{\theta} G]$ is isomorphic to the bundle $\bigoplus_j \Ocal(-d_j)$ over $\Gr(k,{n})$. 

As in the toric case, the critical locus in this phase is given by $p_1 = \cdots = p_s = 0 = F_1 = \cdots = F_s$, so we recover the complete intersection $F_1 = \cdots = F_s$ in $\Gr(k,{n})$, and we call this phase the \emph{geometric phase}.

Just as for the toric complete intersection, we choose the $C^*_R$-action to have weight $0$ on the space $M_{n,k}$ and weight $1$ on all of the $p_j$, so that $W$ has $\C^*_R$ weight $1$ and $\Gamma \cong \GL(k)\times \C^*$.  

The trivial lift $\lift_0$ is a good lift because each monomial of the form $B_{i_1,\dots,i_k}^{\thetaweight} \ell$ is $\Gamma$ invariant for the action induced by $\lift_0$.

The prestable moduli problem of LG-quasimaps to the critical locus $\crst_{\mathrm{geom}}$ consists of maps from prestable orbicurves to the complete intersection
\begin{equation}\label{eq:GR-geom-phase-X}
      \{(\Ccal, f:\Ccal\to Z_{d_1,\dots,d_s})\}.
\end{equation}
If $\Ecal$ denotes the tautological bundle on $\Gr(k,n)$, then the line bundle $\sigma^*(\Lcal_{\lift_0})$ is the $e$th power $\u^*(\Lcal_{\lift_0}) = \det(f^* \Ecal)^e$ of the determinant of the pullback---corresponding to the fact that any $U \in G$ acts on $\ell$ by $\det(U)^{-e}$ and $\C^*_R$ acts on $\ell$ trivially.

\subsubsection{LG-phase}

We call the case where $\thetaweight<0$ the \emph{LG-phase}.  In order to be $G$-invariant, a monomial $\prod B_{i_1,\dots,i_k}^{b_{i_1,\dots,i_k}} \prod p_j^{a_j}\ell^m$ must have $\sum a_j>0$, which implies that any points with every $p_j =0$ must be unstable, but for each $m>0$ and each $j$ the monomial $p_j^{m\thetaweight} \ell^{m d_j}$ is $G$-invariant, so every point with at least one nonzero $p_j$ is $\theta$-semistable.  Therefore $V^{ss}_G(\theta) = M_{k,n} \times (\CC^s\setminus \{\0\})$.  Again, since the $F_j$ and the image of the Pl\"ucker embedding are transverse,  the equations 
$\partial_{B_{i_1,\dots,i_k}}W=\sum_j p_j \partial_{B_{i_1,\dots,i_k}} F_j=0$ imply that the critical locus is
$\left[(\{\0\} \times (\CC^s\setminus \{\0\}))/\GL(k,\CC)\right]$ inside $[V\!\git{\theta} G] = \left[(M_{k,n} \times (\CC^s \setminus \{\0\}))/\GL(k,\CC)\right]$.

This phase does not immediately fit into our theory because we have an infinite stabilizer $\SL(k, \CC)$ for any points of the form $(\0, p_1, \dots , p_s)$. This means that the quotient $[V\!\git{\theta} G]$ is an Artin stack (not Deligne-Mumford).

Hori-Tong \cite{HoTo:07} have analyzed the gauged linear sigma model of the Calabi-Yau complete
intersection $Z_{1, \dots, 1}\subset \Gr(2,7)$ which is defined by seven
linear equations in the Pl\"ucker coordinates. They gave a physical
derivation that its LG-phase is equivalent to the Gromov-Witten theory of
the so-called \emph{Pfaffian variety}
$$\Pf(\bigwedge^2 \CC^7)=\{A\in \bigwedge^2 \CC^7; A\wedge A\wedge
A=0\}.$$ It is interesting to note that the Pfaffian $\Pf(\bigwedge^2 \CC^7)$ is not a complete intersection. For additional work on this example, see \cite{Rod:00, Kuz:08, HoKo:09, ADS:13}

\subsection{Complete Intersections in a Flag Variety}
Another class of interesting examples is that of complete intersections in partial flag varieties.  The partial flag variety $\Fl(d_1,\cdots, d_k)$
parametrizes the space of partial flags
$$0\subset V_1\subset \cdots V_i\subset\cdots V_k=\CC^n$$ such that
$\dim V_i=d_i.$ The combinatorial structure of the equivariant cohomology
of $\Fl(d_1, \cdots, d_k)$ is a very interesting subject in its own right. 

For our purposes, $\Fl(d_1, \cdots, d_k)$ can be constructed as
a GIT or symplectic quotient  of the vector space 
\[
 \prod_{i=1}^{k-1} M_{d_i,d_{i+1}}
\]
by the group 
\[
G = \prod_{i=1}^{k-1} \GL(d_i,\CC)
\]

The moment map sends the element $(A_1,\dots,A_{k-1})\in \prod_{i=1}^{k-1} M_{i,i+1}$ 
to the element $\frac{1}{2}(A_1\bar{A}_1^T,\dots,A_{k-1}\bar{A}_{k-1}^T) \in \prod_{i=1}^{k-1}{\ufrak}({d}_i)$. 

Let the $\chi_i$ be the character of $\prod_j \GL({d}_j)$ given by the determinant of
$i$th factor.  Each character $\chi_i$ defines a line bundle on the vector space $ M_{d_1,d_2}\times \cdots \times M_{d_{k-1},d_k}$, which descends to a line bundle $\XLB_i$ on $\Fl(d_1,\cdots, {n}_k)$. A hypersurface of multidegree $(\ell_1, \dots, \ell_k)$ is a section of $\bigotimes_j \XLB^{\ell_j}_j$. 

To consider the gauged linear sigma model for the complete intersection $F_1=\cdots =F_s=0$ of such sections, we again consider the vector space 
\[
V = \prod_{i=1}^{k-1} M_{d_i,d_i+1} \times \CC^s,
\]
with coordinates $(p_{1},\dots,p_{s})$ on $\CC^s$ and 
superpotential
\[
W=\sum_{j=1}^s p_{j} F_j.
\]
We define an action of $G$ on $p_i$ by $(g_1,\dots,g_{k-1}) \in G$ acts on $p_i$ as $\prod_{j=1}^{k-1} \det(g_j)^{-\ell_{ij}}$, where $\ell_{ij}$ is the $j$th component of the multidegree degree of $F_i$. 

We may describe the polarization as 
\[
\theta = \prod_{i=1}^{k-1} \det(g_i)^{-\thetaweight_i},
\]
or the moment map as 
\[
\mmap(A_1,\dots,A_{k-1},p_1,\dots,p_s) =\frac{1}{2}(A_1\bar{A}_1^T - \sum_{i=1}^s \ell_{1j}|p_j|^2,\dots,A_{k-1}\bar{A}_{k-1}^T - \sum_{i=1}^s \ell_{k-1,j}|p_j|^2 ).
\] 
This gives a phase structure similar to a complete intersection in a product of projective spaces. 

For example, when $\thetaweight_i>0$ for all $i\in \{1,\dots,k-1\}$ we can choose a compatible $\CC^*_R$ action with weight $1$ on $p_j$ and weight $0$
on each $A_i$, and again the trivial lift $\lift_0$ is a good lift of $\theta$ in this phase.

This example should be easy to generalize to complete intersections in quiver varieties.  It would be very interesting to calculate the details of our theory for these examples.

\section{General comments}
When $G$ is non-Abelian or the R-charge is not integral, $G$ and $\CC^*_R$ interact
in a nontrivial way and the description of moduli space is more complicated.  For more details, we encourage readers
to consult \cite{FJR:15a}.  

An important technique that we have not touched on here is cosection localization, which is the main tool for constructing a virtual cycle for the GLSM.  Starting from the noncompact stack of LG-quasimaps to $[V\git{\theta} G]$, the cosection localization technique enables us
to construct a virtual cycle supported on the compact substack of LG-quasimaps to the critical locus of $W$. 

Finally, we remark that the choice of stability condition in our paper is by no means unique. There are other choices of stability conditions that result in different moduli spaces. Please see \cite{CLL:15, ChKi:15} for examples.

\def\cprime{$'$}


\begin{thebibliography}{{CCFK}14}

\bibitem[ADS13]{ADS:13}
N.~Addington, W.~Donovan, and E.~Segal.
\newblock The Pfaffian-Grassmannian equivalence revisited.
\newblock \emph{ArXiv e-prints}, Jan 2014, arXiv:1401.3661.

\bibitem[AJ03]{AbJa:03}
D.~Abramovich and T.~J.~Jarvis.
\newblock Moduli of twisted spin curves.
\newblock \emph{Proc. Amer. Math. Soc.}, 131(3):685--699 (electronic), 2003.

\bibitem[AV02]{AbVi:02}
D.~Abramovich and A.~Vistoli.
\newblock Compactifying the space of stable maps.
\newblock {\em J. Amer. Math. Soc.}, \textbf{15}(1):27--75 (electronic), 2002.

\bibitem[Alp13]{Alp:13} J.~Alper. 
\newblock Good moduli spaces for Artin Stacks.
\newblock \emph{Annales de l'Institut Fourier}, \textbf{63} (2013), no.~6, 2349--2402.

\bibitem[AP10]{AsPl:09} P.S.~Aspinwall and M.R.~Plesser,  
\newblock Decompactifications and Massless D-Branes in Hybrid Models
\newblock \emph{Journal of High Energy Physics}, {(2010)} no.~078. 

\bibitem[Bor91]{Bor:91}
A.~Borel.
\newblock \emph{Linear algebraic groups}, vol. \textbf{126} of \emph{Graduate Texts in
  Mathematics}.
\newblock Springer-Verlag, New York, second edition, 1991.

\bibitem[CCK14]{CCFK:14}  
D.~Cheong, I.~Ciocan-Fontanine, and B.~Kim,
\newblock Orbifold Quasimap Theory.
\newblock \emph{ArXiv e-prints}, May 2014, ArXiv:math/1405.7160.


\bibitem[CG10]{ClGu:10}
P.~{Clarke} and J.~{Guffin}.
\newblock {On the existence of affine Landau-Ginzburg phases in gauged linear
  sigma models}.
\newblock \emph{ArXiv e-prints}, April 2010, arXiv:math/1004.2937.

\bibitem[CN]{CN}
A.~Chiodo and J.~Nagel. In progress.

\bibitem[CIR12]{CIR:12}
A.~{Chiodo}, H.~{Iritani}, and Y.~{Ruan}.
\newblock {Landau-Ginzburg/Calabi-Yau correspondence, global mirror symmetry
  and Orlov equivalence}.
\newblock \emph{Publ. Math. Inst. Hautes \'Etudes Sci.} \textbf{119} (2014), 127--216.

\bibitem[CJR]{CJR}
E. ~{Clader}, F. ~{Janda} and Y.~{Ruan}.
\newblock{In Progress}

\bibitem[CK10]{CFKi:10}
I.~Ciocan-Fontanine and B.~Kim.
\newblock Moduli stacks of stable toric quasimaps.
\newblock \emph{Adv. Math.}, 225(6):3022--3051, 2010.

\bibitem[CKM11]{CFKM:11}
I.~{Ciocan-Fontanine}, B.~{Kim}, and D.~{Maulik}.
\newblock {Stable quasimaps to GIT quotients}.
\newblock \emph{J. Geom. Phys.} \textbf{75} (2014), 17--47.

\bibitem[CL11]{ChaLi:11}
H.-L. {Chang} and J.~{Li}.
\newblock Gromov-witten invariants of stable maps with fields.
\newblock \emph{International Mathematics Research Notices}, RNR 186,
  arXiv:math/1101.0914.

\bibitem[CLL13]{CLL:13}
H.-L. {Chang}, J.~{Li}, and W.~{Li}.
\newblock{Witten's top Chern class via cosection localization}
\newblock To appear in \emph{Invent. Math.}  \emph{ArXiv e-prints}, Mar 2013, arXiv:1303.7126

\bibitem[CLLL15]{CLL:15}
H.-L. {Chang}, J.~{Li}, W.~{Li}, and C.-C. Liu
\newblock{Mixed-spin-$p$-fields of Fermat quintic polynomials}
\newblock \emph{ArXiv e-prints}, May 2015, arXiv:1505.07532v1.

\bibitem[CheR02]{ChRu:02}
W.~Chen and Y.~Ruan.
\newblock Orbifold {G}romov-{W}itten theory.
\newblock In {\em Orbifolds in mathematics and physics ({M}adison, {WI},
  2001)}, volume 310 of {\em Contemp. Math.}, pages 25--85. Amer. Math. Soc.,
  Providence, RI, 2002.

\bibitem[Chi06]{Chi:06a}
A.~Chiodo.
\newblock The {W}itten top {C}hern class via {$K$}-theory.
\newblock \emph{J.~Algebraic Geom.}, 15(4):681--707, 2006.

\bibitem[ChiR10]{ChiRu:10}
A.~Chiodo and Y.~Ruan.
\newblock Landau-{G}inzburg/{C}alabi-{Y}au correspondence for quintic
  three-folds via symplectic transformations.
\newblock \emph{Invent. Math.}, 182(1):117--165, 2010.

\bibitem[ChiR11]{ChiRu:11}
A.~Chiodo and Y.~Ruan.
\newblock L{G}/{CY} correspondence: the state space isomorphism.
\newblock \emph{Adv. Math.}, 227(6):2157--2188, 2011.

\bibitem[{Cla}13]{Cla:13}
E.~{Clader}.
\newblock {Landau-Ginzburg/Calabi-Yau correspondence for the complete
  intersections X\_$\{$3,3$\}$ and X\_$\{$2,2,2,2$\}$}.
\newblock \emph{ArXiv e-prints}, January 2013, arXiv:math/1301.5530.

\bibitem[{Cla}14]{ClaDis:14}
E.~{Clader}.
\newblock {The Landau-Ginzburg/Calabi-Yau Correspondence for
Certain Complete Intersections}.
Doctoral Dissertation, University of Michigan, 2014 



\bibitem[DH98]{DoHu:98}
I.~V. Dolgachev and Y.~Hu.
\newblock Variation of geometric invariant theory quotients.
\newblock \emph{Inst. Hautes \'Etudes Sci. Publ. Math.}, (87):5--56, 1998.
\newblock With an appendix by Nicolas Ressayre.

\bibitem[EJK10]{EJK:10}
D.~Edidin, T.~J.~Jarvis, and T.~Kimura, \emph{Logarithmic trace and
  orbifold products}, Duke Math. J. \textbf{153} (2010), no.~3, 427--473.
  

\bibitem[FJR07]{FJR:07b}
H.~Fan, T.~J.~Jarvis, and Y.~Ruan.
\newblock The {W}itten equation and its virtual fundamental cycle.
\newblock \emph{ArXiv e-prints}, Dec 2007, arXiv:math/0712.4025

\bibitem[FJR08]{FJR:08}
H.~Fan, T.~J.~Jarvis, and Y.~Ruan.
\newblock Geometry and analysis of spin equations.
\newblock \emph{Comm. Pure Appl. Math.}, 61(6):745--788, 2008.

\bibitem[FJR12]{FJR:07a}
H.~Fan, T.~J.~Jarvis, and Y.~Ruan.
\newblock {T}he {W}itten equation, mirror symmetry and quantum singularity
  theory.
\newblock \emph{Annals of Mathematics}
 \textbf{178} (2013), 1--106. %{http://dx.doi.org/10.4007/annals.2013.178.1.1}

\bibitem[FJR15a]{FJR:15a}
H.~Fan, T.~J.~Jarvis, and Y.~Ruan.
\newblock A mathematical theory of the Gauged Linear Sigma Model.
\newblock \emph{Preprint}, 2015. \newblock \emph{ArXiv e-prints}, arXiv:math/1506.02109


\bibitem[FJR15b]{FJR:15b}
H.~Fan, T.~J.~Jarvis, and Y.~Ruan.
\newblock Analytic Theory of the Gauged Linear Sigma Model.
\newblock \emph{Preprint}, 2016. 

\bibitem[FSZ10]{FSZ}
C.~Faber, S.~Shadrin, and D.~Zvonkine.
\newblock Tautological relations and the $r$-spin Witten conjecture.
\newblock \emph{Ann. Sci. \'Ec.~Norm.~Sup\'er.} (4) \textbf{43} (2010), no.~4, 621--658.

\bibitem[CK15]{ChKi:15}
J.~Choi and Y.-H.~Kiem.
\newblock Landau-Ginzburg/Calabi-Yau correspondence via quasi-maps
\newblock Preprint, 2015

\bibitem[GW08]{GoWo:08}
E.~{Gonzalez} and C.~{Woodward}.
\newblock {Quantum Witten localization and abelianization for QDE solutions}.
\newblock November 2008, arXiv:math/0811.3358.


\bibitem[HK09]{HoKo:09}
S.~Hosono and Y.~Konishi.
\newblock Higher genus {G}romov-{W}itten invariants of the {G}rassmannian, and
  the {P}faffian {C}alabi-{Y}au 3-folds.
\newblock \emph{Adv. Theor. Math. Phys.}, 13(2):463--495, 2009.

\bibitem[HK10]{HeKu:10}
C.~{Hertling} and R.~{Kurbel}.
\newblock {On the classification of quasihomogeneous singularities}.
\newblock \emph{J. Singul.} \textbf{4} (2012), 131--153.

\bibitem[HM98]{HaMo:98}
J.~Harris and I.~Morrison.
\newblock {\em Moduli of curves}, volume 187 of {\em Graduate Texts in
  Mathematics}.
\newblock Springer-Verlag, New York, 1998.



\bibitem[HN01]{HoNa:01}
Y.~I. Holla and M.~S. Narasimhan.
\newblock A generalisation of {N}agata's theorem on ruled surfaces.
\newblock \emph{Compositio Math.}, 127(3):321--332, 2001.

\bibitem[HT07]{HoTo:07}
K.~Hori and D.~Tong.
\newblock Aspects of non-abelian gauge dynamics in two-dimensional
  {$\mathscr{N}=(2,2)$} theories.
\newblock \emph{J.~High Energy Phys.}, (5):079, 41 pp. (electronic), 2007.

\bibitem[JKV01]{JKV} T. Jarvis, T. Kimura and A. Vaintrob, \emph{Moduli spaces
of higher spin curves and integrable hierarchies}, Compositio
Math. \textbf{126} (2001), 157--212.


\bibitem[KL13]{KiLi:10}
Y.-H. {Kiem} and J.~{Li}.
\newblock {Localizing Virtual Cycles by Cosections}.
\newblock \emph{Jour. Amer. Math. Soc.} \textbf{26} (2013), No. 4, 1025--1050.
%\newblock \emph{ArXiv e-prints}, July 2010, arXiv:math/1007.3085.

\bibitem[KS11]{KrSh:11}
M.~{Krawitz} and Y.~{Shen}.
\newblock {Landau-Ginzburg/Calabi-Yau Correspondence of all Genera for Elliptic
  Orbifold {$\mathbb{P}^1$}}.
\newblock June 2011, arXiv:math/1106.6270.

\bibitem[Kau06]{Kau:06}
R.~M. Kaufmann.
\newblock Singularities with symmetries, orbifold {F}robenius algebras and
  mirror symmetry.
\newblock In \emph{Gromov-{W}itten theory of spin curves and orbifolds}, volume
  403 of \emph{Contemp. Math.}, pages 67--116. Amer. Math. Soc., Providence, RI,
  2006.

\bibitem[{Kim}11]{Kim:11}
B.~{Kim}.
\newblock {Stable quasimaps}.
\newblock \emph{ArXiv e-prints}, June 2011, arXiv:math/1106.0804.

\bibitem[Kon92]{Kon:92}
M.~Kontsevich.
\newblock Intersection theory on the moduli space of curves and the matrix
  {A}iry function.
\newblock \emph{Comm. Math. Phys.}, 147(1):1--23, 1992.

\bibitem[Kra10]{Kra:10}
M.~Krawitz.
\newblock \emph{{F}{J}{R}{W} Rings and {L}andau-{G}inzburg Mirror Symmetry}.
\newblock PhD thesis, University of Michigan, 2010.

\bibitem[Kuz08]{Kuz:08}
A.~Kuznetsov.
\newblock Lefschetz decompositions and categorical resolutions of
  singularities.
\newblock \emph{Selecta Math. (N.S.)}, 13(4):661--696, 2008.

\bibitem[LPS]{LPS}
Y.~{Lee}, N.~{Priddis} and M.~{Shoemaker}.
\newblock A proof of the Landau-Ginzburg/Calabi-Yau correspondence via the crepant transformation conjecture.
\newblock \emph{ArXiv e-prints}, October 2014,  arXiv:1410.5503

\bibitem[Lee06]{Lee}
Y.-P.~Lee.
\newblock {W}itten's conjecture and the {V}irasoro conjecture for genus up to
  two.
\newblock In {\em Gromov-{W}itten theory of spin curves and orbifolds}, volume
  403 of {\em Contemp. Math.}, pages 31--42. Amer. Math. Soc., Providence, RI,
  2006.


\bibitem[MFK94]{MFK:94}
D.~Mumford, J.~Fogarty, and F.~Kirwan.
\newblock \emph{Geometric invariant theory}.
\newblock Springer-Verlag, Berlin, third edition, 1994.

\bibitem[MOP11]{MOP:11}
A.~Marian, D.~Oprea, and R.~Pandharipande.
\newblock The moduli space of stable quotients.
\newblock \emph{Geom. Topol.}, 15(3):1651--1706, 2011.


\bibitem[Ols07]{Ols:07}
M.~C. Olsson.
\newblock ({L}og) twisted curves.
\newblock {\em Compos. Math.}, 143(2):476--494, 2007.

\bibitem[PS]{PS}
N.~{Priddis} and M.~{Shoemaker}
\newblock{A Landau-Ginzburg/Calabi-Yau correspondence for the mirror quintic}.
\newblock \emph{ArXiv e-prints}, September 2013, arXiv:1309.6262 

\bibitem[PV01]{PoVa:01}
A.~Polishchuk and A.~Vaintrob.
\newblock Algebraic construction of {W}itten's top {C}hern class.
\newblock In \emph{Advances in algebraic geometry motivated by physics
  ({L}owell, {MA}, 2000)}, volume 276 of \emph{Contemp. Math.}, pages 229--249.
  Amer. Math. Soc., Providence, RI, 2001.

\bibitem[PV11]{PoVa:11}
A.~{Polishchuk} and A.~{Vaintrob}.
\newblock {Matrix factorizations and Cohomological Field Theories}.
\newblock \emph{J. Reine Angew. Math.} (2014) publ. online DOI: 10.1515/crelle-2014-0024, 128 pages.
%{ArXiv e-prints}, May 2011, arXiv:math/1105.2903.

\bibitem[{Pop}12]{Pop:12}
A.~{Popa}.
\newblock {Two-Point Gromov-Witten Formulas for Symplectic Toric Manifolds}.
\newblock \emph{ArXiv e-prints}, June 2012, arXiv:1206.2703.

\bibitem[RR]{RR}
D.~{Ross} and Y.~{Ruan}
\newblock{Wall-Crossing in Genus Zero Landau-Ginzburg Theory}
\newblock \emph{ArXiv e-prints}, February 2014, arXiv:1402.6688 

\bibitem[RRS]{RRS}
D.~{Ross}, Y.~{Ruan} and M. ~{Shoemaker}
\newblock{In Progress}

\bibitem[Ram96]{Ram:96a}
A.~Ramanathan.
\newblock Moduli for principal bundles over algebraic curves. {I}.
\newblock \emph{Proc. Indian Acad. Sci. Math. Sci.}, 106(3):301--328, 1996.

\bibitem[R{\o}d00]{Rod:00}
E.~A. R{\o}dland.
\newblock The {P}faffian {C}alabi-{Y}au, its mirror, and their link to the
  {G}rassmannian {$G(2,7)$}.
\newblock \emph{Compositio Math.}, 122(2):135--149, 2000.


\bibitem[Rua12]{Rua:12}
Y.~Ruan.
\newblock The {W}itten equation and the geometry of the {L}andau-{G}inzburg
  model.
\newblock \emph{Preprint}, June 2012.

\bibitem[TX]{TX}
G. ~{Tian} and G. ~{Xu}
\newblock{Analysis of gauged Witten equation}.
\newblock  arXiv:1405.6352

\bibitem[Tay11]{ComDiag:11}
P.~Taylor.
\newblock Commutative diagrams package.
\newblock {http://http://www.paultaylor.eu/diagrams/}, 2011.

\bibitem[Tha96]{Tha:96}
M.~Thaddeus.
\newblock Geometric invariant theory and flips.
\newblock \emph{J.~Amer. Math. Soc.}, 9(3):691--723, 1996.

\bibitem[VW89]{VaWa:89}
C.~{Vafa} and N.~{Warner}.
\newblock {Catastrophes and the classification of conformal theories}.
\newblock \emph{Physics Letters B}, 218:51--58, February 1989.

\bibitem[Wit91]{Wit:91}
E.~Witten.
\newblock Two-dimensional gravity and intersection theory on moduli space.
\newblock In \emph{Surveys in differential geometry ({C}ambridge, {MA}, 1990)},
  pages 243--310. Lehigh Univ., Bethlehem, PA, 1991.

\bibitem[{Wit}92]{Wit:92a}
E.~{Witten}.
\newblock {The {$N$}-matrix model and gauged {WZW} models}.
\newblock \emph{Nuclear Physics B}, 371:191--245, March 1992.

\bibitem[Wit93]{Wit:93}
E.~Witten.
\newblock Algebraic geometry associated with matrix models of two-dimensional
  gravity.
\newblock In \emph{Topological methods in modern mathematics ({S}tony {B}rook,
  {NY}, 1991)}, pages 235--269. Publish or Perish, Houston, TX, 1993.

\bibitem[Wit97]{Wit:97}
E.~Witten.
\newblock Phases of {$N=2$} theories in two dimensions.
\newblock In \emph{Mirror symmetry, {II}}, volume~1 of \emph{AMS/IP Stud. Adv.
  Math.}, pages 143--211. Amer. Math. Soc., Providence, RI, 1997.

\end{thebibliography}
\end{document}